\documentclass[onefignum,onetabnum]{siamonline220329}
\usepackage[utf8]{inputenc}
\usepackage[english]{babel}
\usepackage{bbm}
\usepackage{amsfonts}
\usepackage{amssymb}
\usepackage{graphicx}
\usepackage{booktabs}
\usepackage{subcaption}
\usepackage{multirow}
\usepackage[algo2e,ruled,vlined]{algorithm2e}
\usepackage{enumitem}
\setlist[enumerate]{leftmargin=.5in}
\setlist[itemize]{leftmargin=.5in}
\usepackage{float}
\usepackage{dsfont}

\def\R {{\mathbb R}}

\def\N {{\mathbb N}}

\def\X {{\mathbb X}}

\def\Bu {{\mathbf{u}}}
\def\Bv {{\mathbf{v}}}
\def\Bw {{\mathbf{w}}}
\def\Bx {{\mathbf{x}}}
\def\By {{\mathbf{y}}}

\def\Bq {{\mathbf{q}}}
\def\Bp {{\mathbf{p}}}

\def\cB {{\cal B}}

\def\cD {{\cal D}}

\def\cP {{\cal P}}
\def\cR {{\cal R}}

\def\cF {{\cal F}}
\def\cH {{\cal H}}

\def\cL {{\cal L}}
\def\cH {{\cal H}}
\def\cN {{\cal N}}
\def\cM {{\cal M}}
\def\cS {{\cal S}}
\def\cT {{\cal T}}
\def\cW {{\cal W}}
\def\cX {{\cal X}}
\def\cY {{\cal Y}}

\def\Pp {\Bp_{l+1}}
\def\Pl {\Bp_{l}}
\def\Xp {\Bx_{l+1}}
\def\Xl {\Bx_l}

\DeclareMathOperator{\co}{conv}

\def\eps {\epsilon}

\newcommand{\normiii}[1]{{\vert\kern-0.25ex\vert\kern-0.25ex\vert #1 
    \vert\kern-0.25ex\vert\kern-0.25ex\vert^2}}

\makeatletter
\newcommand{\mylabel}[2]{#2\def\@currentlabel{#2}\label{#1}}

\DeclareMathOperator*{\argmin}{arg\,min}

\DeclareMathOperator*{\argmax}{arg\,max}

\DeclareMathOperator*{\conv}{conv}

\DeclareMathOperator*{\sgn}{sgn}

\newsiamremark{remark}{Remark}

\begin{document}
\title{\bf The Pontryagin Maximum Principle for Training Convolutional Neural Networks}
\author{S. Hofmann 
\thanks{Institut f\"ur Mathematik, Universit\"at W\"urzburg, Emil-Fischer-Strasse 30,
97074 W\"urzburg, Germany. \email{ 
sebastian.hofmann@uni-wuerzburg.de}}
\and A. Borz{\`i}
\thanks{Institut f\"ur Mathematik, Universit\"at W\"urzburg, Emil-Fischer-Strasse 30,
97074 W\"urzburg, Germany. \email{ alfio.borzi@uni-wuerzburg.de}}}

\headers{The Pontryagin Maximum Principle for Training CNNs}{S. Hofmann, and A. Borz{\`i}}

\maketitle

\begin{abstract}
 
A novel batch sequential quadratic Hamiltonian (bSQH) algorithm for training convolutional neural networks (CNNs) with $L^0$-based regularization is presented. This methodology is based on a discrete-time Pontryagin maximum principle (PMP). It uses forward and backward sweeps together with the layerwise approximate maximization of an augmented Hamiltonian function, where the augmentation parameter is chosen adaptively.

A technique for determining this augmentation parameter is proposed, and the loss-reduction and convergence properties of the bSQH algorithm are analysed theoretically and validated numerically.  
Results of numerical experiments in the context of image classification with a sparsity-enforcing $L^0$-based regularizer demonstrate the effectiveness of the proposed method in full-batch and mini-batch modes.

 \end{abstract}
 


 \bigskip
 	\section{Introduction}
 	\label{sec:intro}
In the last decades, convolutional neural networks (CNNs) have become a standard choice for computer vision tasks and, in particular, medical imaging analysis, where they are applied to aid in detecting and diagnosing diseases and classifying tumour types based on MRI, CT  or X-ray scans.  Compared to its predecessors, like the multilayer perceptron (MLP), CNNs utilize convolutional layers with trainable filters in addition to fully connected layers, making them especially suited for processing data with a grid-like topology. This network architecture enables a more efficient extraction of local patterns and relevant features from an input image, by simultaneously requiring fewer parameters than fully connected architectures. The development of CNNs started in the 1980s when Kunihiko Fukushima \cite{Fukushima1980} introduced the 'Neocognitron' as the earliest precursor of convolutional networks. Almost one decade later, Yann LeCun \cite{LeCun1988} developed the LeNet-5 architecture for document recognition, representing a milestone in the design of modern CNNs. Over the last decades, a variety of architectures have contributed to the success of these network types, amongst which there is AlexNet \cite{AlexNet2012}, which popularized the use of ReLU activation functions and VGGNet \cite{VGGNet2015} where an increased network depth was used, showing the benefits of using smaller filter sizes. However, it occurs that very deep architectures often face a vanishing gradient problem, preventing very deep networks from being trained efficiently using gradient-based training methods and backpropagation. To mitigate this vanishing gradient problem, Kaiming He et al. \cite{ResNetProp} developed the ResNet, a network architecture that uses skip connections and thereby enables the training of even deeper neural networks. 

Apart from this, the introduction of ResNet-like architectures gave birth to the dynamical system viewpoint of deep learning, where networks with a residual structure are related to numerical schemes for solving differential equations; see \cite{EWeinan2017,Haber2017,Chang2017,Chen2018,Lu2018,Ruthotto2020}. The dynamical systems viewpoint offers important insights into the stability of a network's forward propagation process by allowing the study of very deep networks from a continuous-time dynamical systems perspective. More importantly for our investigation, it motivates the development of novel learning techniques based on optimal control theory \cite{EWeinanEMSA2018,Hofmann2022}.

While serving as a motivation for developing new training algorithms, the continuous-time dynamical systems interpretation can be limiting, particularly when considering network architectures without skip connections, such as LeNet-5 and AlexNet. In those cases, the continuous-time interpretation becomes even more difficult when taking into account that CNN models use dimensionality reduction as a design approach for feature extraction. However, even if the relation to the continuous-time setting is lost for most modern network architectures, it has been shown that special optimal control-based techniques can still be applied to the machine learning framework, almost independently of the network's structure \cite{LiHaoEMSA2018}. Based thereon, in this work, we considerably extend our optimal control-based learning technique presented in \cite{Hofmann2022} to the framework of CNNs. For this purpose, we construct a CNN as a discrete-time controlled system including trainable filters, dimensionality reduction in the form of pooling operations, and a sequence of fully-connected layers for classification. Further, for this extension and within this special structure, we also take into account the computational complexity arising from the increased number of parameters and dimensionality of data, which can impose memory constraints. To address these challenges, we derive a novel mini-batch variant of our optimal control-based learning algorithm. Furthermore, we include sparse network training in the discussion of our algorithm, where the aim is to determine the essential parameters for a given model by forcing irrelevant or redundant ones to be zero.  

In the context of CNNs, where models are often over-parametrized, deriving sparse solutions can significantly reduce the memory usage of a trained model while improving the computational efficiency of the forward propagation process. In the classical backpropagation framework, commonly a $L^1$-norm penalty is considered to induce sparsity. However, we demonstrate that with our optimal control-based technique, also $L^0$-'norm' regularized machine learning problems can be solved. In this way, we can determine CNN models with a higher level of sparsity in the trained parameters compared to the $L^1$-based Elastic-Net approach \cite{Zou2005}, while simultaneously achieving a similar level of accuracy.

In the next section, we introduce the interpretation of supervised learning tasks with CNN constraints as a discrete-time optimal control problem, which enables us to formulate first-order optimality conditions in the form of a discrete Pontryagin maximum principle (PMP). The discrete PMP, as formulated by Halkin \cite{Halkin1966}, is an extension of the work by Pontryagin et al. \cite{PBG1962}, who proved the PMP in the continuous-time optimal control framework; see also \cite{BGP1956}.
Opposed to the classical machine learning setting, where the first-order optimality conditions are based on the conditions of Karush-Kuhn-Tucker (KKT); see, for example, \cite{Aggarwal2018,NocedalWright2006}, the PMP states that optimal trainable parameters are characterized by maximizing an auxiliary function at each layer, which is called the Hamilton-Pontryagin (HP) function. This maximization property allows the formulation of optimality conditions in the supervised learning framework even without differentiability assumptions on the loss function with respect to the trainable parameters. In Section \ref{sec-Preliminaries and Assumptions}, we formulate and discuss requirements for our PMP-based learning approach in the classical supervised learning framework, which is followed by a discussion of the discrete PMP for machine learning with CNNs in Section \ref{sec-The maximum principle}. In Section \ref{sec-Successive approximation schemes}, we reason that the PMP motivates the formulation of a generalized backpropagation algorithm for CNNs based on a discrete version of the method of successive approximations (MSA). This early method was initially introduced by Krylov and Chernous'ko \cite{Krylov1963} for continuous-time optimal control problems, see also \cite{Kelley1961}, and has since then undergone further improvements and extensions \cite{EWeinanEMSA2018,Hofmann2022,LiHaoEMSA2018,SakawaShindo1980,ChernouskoLyubushin1982,BreitenbachBorzi2018}.
In contrast to the original MSA method, these more recent modifications use an augmented HP function with a weighted penalization term to accomplish a stable optimization process. In particular, the choice of this augmentation weight has a decisive influence on the convergence of the algorithms, which makes methods such as the sequential quadratic Hamiltonian (SQH) scheme \cite{Hofmann2022} particularly suitable for the definition of a stable training procedure, since it automatically determines an augmentation weight which guarantees a monotonic decrease in the loss functional value. We present a novel batch SQH method (bSQH), which enables efficient training of CNN architectures in a mini-batch fashion and improves upon the process of estimating the required augmentation weight incorporated into the augmented HP function. A decisive advantage of our methodology compared to the standard backpropagation approach lies in the formulation of the updating step, which does not require the existence of gradients with respect to the trainable parameters. Consequently, our framework can accommodate supervised learning even with non-continuous regularizers and also enables the training of neural networks where the set of admissible parameters is discrete \cite{LiHaoEMSA2018,Baldassi2007,Baldassi2017}.

We present in Section \ref{sec-Sparse network training} a reformulation of non-residual CNN architectures such as LeNet-5 that enables the derivation of an analytic expression for the trainable parameters maximizing the augmented HP functions at each iteration. This allows us to establish the bSQH scheme as a suitable algorithm for training sparse CNNs with a $L^0$-'norm' based regularization technique. In Section \ref{Sec-Well-definedness and Convergence}, we prove the minimizing property of the bSQH algorithm by introducing a sufficient decrease condition for the mini-batch loss. Moreover, in addition to considerably extending our initial results from \cite{Hofmann2022} to more general network structures and sparse learning techniques, we provide a novel convergence result. In particular, this result guarantees convergence of the bSQH scheme to trainable parameters satisfying the necessary optimality conditions of the discrete PMP. We validate the feasibility of the necessary conditions for convergence in Section \ref{sec-numerical experiments}. Further, in this section, we present numerical experiments of sparse network training in the framework of image classification with CNNs and an application to medical imaging. We numerically investigate a novel estimation strategy for the augmentation weights of the augmented HP function, which leads to a considerable reduction in the number of necessary inference steps required for the training process. In this setting, we also compare the bSQH method in its full-batch setting with the mini-batch variant and discuss the computational advantages of the mini-batch approach. To conclude our numerical investigation, our PMP-based training procedure is successfully applied to solve a $L^0$-regularized supervised learning problem, where the involved CNN uses ReLU activations, batch normalization and max-pooling layers to classify abdominal CT images.   
A section of conclusion completes this work.

	\section{The supervised learning problem for CNNs}	
	\label{sec-SVL as OPC}

\begin{figure}[H]
\begin{center}
\includegraphics[width=10cm]{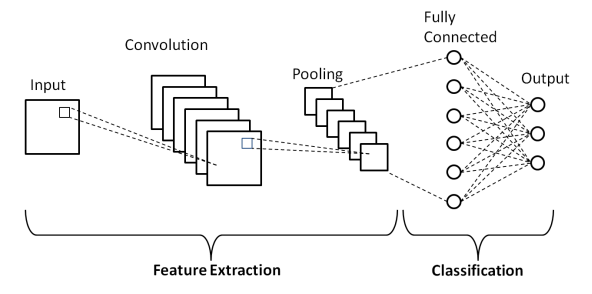}    
\end{center}
\caption{Convolutional neural network architecture \cite{ImageConv}.}
\label{Fig: CNN}
\end{figure}

In the standard supervised learning setting, we have a dataset composed of data-label pairs 
$
\cD: = \mathcal{X} \times \mathcal{Y}, 
$ 
where $\mathcal{X} \subset \R^{n \times c}$ is a bounded set of input data and $\mathcal{Y} \subset \R^m$ is the corresponding set of labels. When considering the typical image classification task, each $x^s \in \mathcal{X}$ is interpreted as a suitably preprocessed and vectorized image, where the dimension $n$ is the total number of pixels in each of the $c$ channels. A standard coloured image, for example, has three channels, each dedicated to one of the primary colours, red, green and blue. For our theoretical investigation, we can assume without loss of generality that $c = 1$ holds by supposing a suitable vectorization. Each $y^s\in \mathcal{Y}$ is a one-hot encoded label representing the affiliation of the image $x^s$ to one of the $m$ different classes of images contained in $\mathcal{X}$. We suppose that there exists an unknown, well-defined and continuous mapping $\mathcal{M}: \mathcal{X} \rightarrow \mathcal{Y}$ describing the data-label relation of elements in $\cD$.
Then, based on the information encoded in the dataset, the supervised learning task is to find a suitable set of parameters $ \Bu=(u_l)_{l=0}^{L-1} $, where each $u_l$ is the element of an appropriate closed Euclidean subset $U_l$, such that the following parametric, composite function approximates the mapping $\mathcal{M}$ up to a satisfying degree
\begin{equation}\label{Function composition}
\cN_{\Bu}: \cX \to \cY, \quad x \mapsto \cN_{\Bu}(x) = \lambda_{L-1}(\;\cdot\;,u_{L-1}) \circ \lambda_{L-2}(\;\cdot\;,u_{L-2}) \circ \ldots \circ \lambda_{0}(x,u_0).
\end{equation}
The composite structure \eqref{Function composition} itself is called the neural network and consists of $L\in \N$ layers $\lambda_l: \R^{n_l} \times U_l \rightarrow \R^{n_{l+1}}$ which, in the framework of CNNs, take the following form
\begin{equation}\label{NN Layer}
  \lambda_l: \R^{n_l} \times U_l \rightarrow \R^{n_{l+1}}\quad \lambda_l(x,u) = \cP_l( \cS_l(x) + \sigma_l (\cW_l(u) x + \cB_l u)).  
\end{equation}
In \eqref{NN Layer}, the transformation $\cS_l: \R^{n_{l}} \rightarrow \R^{\Tilde{n}_{l}}$  defines a projection or identity mapping for modeling skip connections and $\cP_l:\R^{\Tilde{n}_{l}} \rightarrow \R^{n_{l+1}}$ describes the process of a pooling operation or batch normalization \cite{Aggarwal2018}. Pooling operations resemble the application of non-trainable filters, such as maximum or average filters and serve the purpose of reducing the spatial dimension of the processed data while enhancing the network's ability to extract important features. Based upon the universal approximation theorem of George Cybenko \cite{Cybenko1989}, the network's core is represented by the non-linearity $\sigma_l (\cW_l(u) x + \cB_l u)$,  where for each layer $l = 0,\ldots,L-1$, the transformation $\sigma_l: \R^{\Tilde{n}_l} \rightarrow \R^{\Tilde{n}_l}$, is a non-linear component-wise defined activation function. To model the linear affine transformation of the states in \eqref{NN Layer}, we introduce the linear transformations $\cW_l: U_l \rightarrow \R^{\Tilde{n}_l\times n_l}$ and $\cB_l: U_l \rightarrow \R^{\Tilde{n}_l}$ which determine the type of a layer. A CNN is commonly constructed with two types of layers: Convolutional layers for feature extraction and fully connected layers for classification of the extracted features; see Figure \ref{Fig: CNN}. By taking into account that discrete convolutions can be reformulated as matrix multiplication \cite{Aggarwal2018}, in each convolutional layer, the transformation $\cW_l(u)$ takes the form of a sparse Toeplitz matrix. In contrast, for fully connected layers, this transformation is represented by a dense connectivity matrix, where each element acts as an independently adjustable weight.
To determine the optimal parameters $\Bu$ of the architecture \eqref{Function composition} for a given task, a suitable optimization problem is formulated. This so-called supervised learning problem is based on the assumption that, for a given dataset, an appropriate probability distribution $\mu$ on $\cD$ exists such that the abstract supervised learning task can be formulated as follows:
\begin{equation*}\label{General optimization problem}
\begin{aligned}
    \min_{\Bu \in U_{ad}} J(\Bu) &:= \mathbb{E}_{x \sim \mu}[\phi(\cM(x),\cN_{\Bu}(x))] + \rho \;\cR(\Bu),
\end{aligned}
\end{equation*}
where $\cM(x)$ as well as $\cN_{\Bu}(x)$ are functions of random variables with the probability distribution $\mu$. The set of admissible parameters $\Bu$ is defined by 
\begin{equation}\label{U_ad^L}
U_{ad} := U_0 \times U_1 \times \; \cdots \; \times U_{L-1}.
\end{equation}
Further $\cR: U_{ad} \rightarrow \R$ defines a regularizer which is weighted using the fixed scaling factor $\rho \geq 0$ and is applied to improve the stability of the model \eqref{Function composition} or to introduce additional features such as sparsity into the trained parameters. Further, a suitable loss function $\phi: \cY \times \cY \rightarrow \R$ is defined, which attains its minimum if a perfect fit is achieved, i.e. $\cM(x) = \cN_{\Bu}(x)$ for all $x \in \cX$. In practice, only a finite number of observations is available in the form of a data-batch $B \subset \cD$, so that the true probability distribution $\mu$ is unknown and only a finite set of evaluations of $\cM$ exists. Under the assumption that every datapoint in $B$ is independent and identically distributed (i.i.d), the law of large numbers motivates approximating $\mu$ using an empirical measure, which leads to the so-called empirical risk minimization problem:
\begin{equation}\label{Empirical Risk minimization}
\min_{\Bu \in U_{ad}} J_B(\Bu):= \cL_B(\Bu) + \rho \; \mathcal{R}(\Bu),
\end{equation}
where $\cL_B(\Bu) := \frac{1}{|B|}\sum_{s\in I_B} (\Phi_s \circ \cN_{\Bu})(x^{s})$. For each $s \in I_B$, we define 
$
 \Phi_s: \cY \rightarrow \R,\; y \mapsto \phi(y^s,y),
$
where $y^s = \cM(x^s)\in \cY$ denotes the label corresponding to a specific sample $x^s \in B$ with corresponding index $s \in I_B$.

\section{Preliminaries and Assumptions}
	\label{sec-Preliminaries and Assumptions} 

For our discussion, we denote the standard Euclidean norm by $\|\cdot\|$ and the corresponding induced matrix norm with $\|\cdot\|_M$. We further introduce for each $\Bu \in U_{ad}$ the notation $\normiii{\Bu}: = \sum_{l=1}^{L-1} \|u_l\|^2$. To avoid notational clutter, we will henceforth denote any generic constant with $C>0$.
The reachable set of a neural network at layer $l = 0,\ldots, L-1 $, with respect to a given batch of data $B$ is defined as follows 
$$
\X^B_{l} := \bigotimes_{s\in I_B}\X^s_{l}, \qquad \X^s_{l} := \co\left(\{ y \in \R^{n_{l}} \; | \; \exists \; v \in U_{l-1} \; s.t. \; y = \lambda_{l-1}(x,v), \; x \in \X^s_{l-1} \}\right),
$$
where $\X^s_{0} = \{x^s\}$ and '$\conv$' denotes the convex hull operator. Correspondingly, to emphasize the dependency of a feature map upon a fixed set of trainable parameters, we denote with $\Bx_l(\Bu) \in \X_l^{B}$ the entirety of feature maps at a given layer, generated using the trainable parameters $\Bu \in U_{ad}$. Based thereon, the following compact notation is introduced for each layer function corresponding to a given data batch 
$$
\Lambda_l(\Bx_l(\Bu),u_l) := \begin{bmatrix} \lambda_{l}(x^{s_1}_{l}(\Bu),u_l) & \cdots & \lambda_l(x^{s_{|B|}}_{l}(\Bu),u_l) \end{bmatrix}^T,
$$
together with the batch dependent loss function $\Phi_B(\Bx_{L}(\Bu)) := \sum_{s \in I_B} \Phi_s( x^s_{L}(\Bu)) $. To invoke subdifferential calculus, we define, corresponding to \cite{Clarke2013}, for a generic locally Lipschitz function $f:U_l \rightarrow \R^n$ its generalized directional derivative at $\Bu$ in direction $\Bv$ as follows
\begin{equation}\label{generalized directional derivative}
\mathring{f}(u; v) := \limsup_{w \rightarrow u, t \rightarrow 0^+} \frac{f(w + tv) - f(w)}{t}, 
\end{equation}
and denote the corresponding Clarke subdifferential at $u \in U_l$ by
\begin{equation}\label{Ch 4.1: Clarke Subdifferential}
\partial^C f(u) := \{w \in U_l\;|\;  \mathring{f}(u;v) \geq \langle w, v \rangle \quad \forall v \in U_l  \}. 
\end{equation}
For our PMP-based machine learning approach, we make the following Lipschitz assumptions:
\begin{enumerate}
\item[\mylabel{A1}{\textbf{A1}}]
The function $\Phi_B \in C^{1}(\X_L^B,\R)$ is bounded from below and there exists a constant $C>0$ such that for all $\Bx,\By \in \X_L^B$ it holds 
$$
\| \nabla\Phi_B(\Bx) - \nabla\Phi_B(\By)\| \leq C\| \Bx - \By \|.
$$
\item[\mylabel{A2}{\textbf{A2}}]
For any $ l = 0\ldots,L-1$ and all $u_l \in U_l$, we have $\Lambda_l(\cdot, u_l) \in C^1(\X_l^B,\R)$. Further, there exists a constant $C>0$ such that for all $\Bx,\By \in \X_l^B$ and $u,v \in U_l$ it holds
$$
\| \Lambda_l(\Bx,u) - \Lambda_l(\By,v) \| + \| \partial_x\Lambda_l(\Bx,u) - \partial_x\Lambda_l(\By,v)\|_M \leq C(\| \Bx - \By \| + \| u - v \|),
$$ 
for all $l = 0\ldots,L-1$.
\end{enumerate}
The regularity assumption in \ref{A1} on the terminal loss function is satisfied by the state-of-the-art mean squared loss used for regression tasks and the cross-entropy loss commonly applied in classification with CNNs \cite{Aggarwal2018}. To consider the applicability of the Lipschitz assumptions in \ref{A1}-\ref{A2} to the standard neural network setting, we prove the following Proposition for non-residual structures.
\begin{proposition}\label{Proposition-boundedness of x}
For all $l = 0,\ldots,L-1$ in \eqref{NN Layer}, let $\cS_l \equiv 0$ and suppose that $f_l:=\cP_l \circ \sigma_l \in C^{1,1}(\R^{n_l})$ has a Lipschitz derivative. Then if $U_l$ is compact, the Assumption \ref{A2} holds.
\end{proposition}
\begin{proof}
For ease of notation, we consider the single sample case $|B| = 1$ and merge the bias function into the definition of $\cW_l$ by adjusting the network's states correspondingly. Thus, without loss of generality, we can suppose $\cB \equiv 0$. Using elementary Lipschitz estimates and the boundedness of the Jacobian $\partial f_l$ by assumption, we have that for all $(\Bx,u), (\By,v) \in \X^B_l \times U_l $ there exists a constant $C>0$ such that  
\begin{equation}
\| \partial_x \Lambda_l(\Bx,u) - \partial_x \Lambda_l(\By,v) \|_M 
\leq C \left(\|\cW_l(u)\|^2_M \| \Bx - \By \|  +  \|\cW_l(u_l)\|_M\|\By\|\| u - v\| \right).
\end{equation}
 %
A similar reasoning applies to obtain a Lipschitz estimate for the layer functions $\Lambda_l$. The compactness of $U_l$ together with the Lipschitz property of the mapping $\cW_l$ guarantees the boundedness of $\|\cW_l(u)\|_M$ and it remains to prove the boundedness of each $\X_l^B$. By definition, we have $\Bx_{l+1}(\Bu) = f_l(\cW_l(u_l) \Bx_{l}(\Bu))$ and by the Lipschitz property of $f_l$ for all $l = 0, \ldots, L-1$, there exists a constant $C>0$ such that the following inequality holds: 
$$
\| \Bx_{l+1}(\Bu)\| \leq C \|\Bx_l(\Bu)\| + \|f_l(0)\|.
$$
The application of a discrete formulation of Gronwall's inequality reveals the existence of a constant $C>0$ such that 
$$
\| \Bx_{l+1}(\Bu)\| \leq C\left( \|x^s\| + \sum_{l=0}^{L-1}\|f_l(0)\|\right),
$$
proving the boundedness of each $\X_l^B$.
\end{proof}

If $\rho>0$ and the regularizer $\cR$ is lower semi-continuous and coercive, Proposition \ref{Proposition-boundedness of x} ensures that Assumptions \ref{A1}-\ref{A2} can be accommodated in our neural network setting. That is, coupled with the boundedness from below of $\Phi_B$, the properties of the regularizer guarantee that the corresponding level sets $\{ \Bu \;|\; J_{B}(\Bu) \leq C\} \subseteq U_{ad}$ are compact for every $C > 0$. Consequently, the existence of a minimizer for \eqref{Empirical Risk minimization} is assured and the Lipschitz assumptions in \ref{A1}-\ref{A2} can be accommodated by selecting $U_{ad}$ as an appropriate level set.

The case where $\rho = 0$ is also well-posed. In fact, the existence of a minimizer for \eqref{Empirical Risk minimization} can be guaranteed through a so-called max-norm regularization technique, which was initially introduced for CNNs in \cite{Srivastava2014}. In our context, this technique can be modelled by defining $\cW_l$ as the composition of a projection operation onto a given compact and convex set, followed by a linear transformation. This adjustment not only preserves the Lipschitz property of the weight function but also ensures the boundedness of $\X^B$, and thus the validity of \ref{A1}-\ref{A2}. Additionally, due to the structure of our neural network, the regularity assumption in \ref{A2} regarding the states is satisfied for classical sigmoidal activation functions and average pooling operations. We note that, even though functions like ReLU and max-pooling operations are not included, these functions can be smoothly approximated up to any desired degree of accuracy by, e.g. Softplus for ReLU and a weighted Softmax function for max-pooling. For a comprehensive overview of commonly used activation functions, we refer to \cite{Dubey2022}. Consequently, we expect that our results will hold at least asymptotically for ReLU-based networks, ensuring the practical application of our method to a variety of cases. This is supported by numerical evidence presented in Section \ref{sec-numerical experiments}.

  	\section{The maximum principle for CNNs}
  	\label{sec-The maximum principle}

The state-of-the-art approach in machine learning for solving the problem \eqref{Empirical Risk minimization} is to apply a so-called backpropagation technique, introduced by Rumelhart et al. \cite{Rumelhart1986}. This scheme relies on gradient-based optimization methods and is, in its essence, a combination of the chain rule and classical gradient descent \cite{Aggarwal2018,NocedalWright2006}. Commonly, for this approach, it is assumed that the empirical risk function $J_B$ is differentiable. Various derivations of the backpropagation methodology have been introduced, leading to multiple variants of the original algorithm. For our framework, the approach by LeCun \cite{LeCun1988} is especially outstanding as it motivates backpropagation through optimal control-based techniques by invoking variational principles and the Lagrange formalism to formulate first-order necessary conditions for supervised learning problems of the form \eqref{Empirical Risk minimization}. Under suitable differentiability assumptions, these first-order necessary conditions are formulated as follows:

Suppose $\Bu^*\in U_{ad}$ solves \eqref{Empirical Risk minimization} with corresponding feature maps $\Bx^*:=(\Bx_l(\Bu^*))_{l = 0}^{L}$. Then there exists a set of adjoint vectors $\Bp^*:=(\Bp_l(\Bu^*))_{l = 0}^{L}$, such that the following set of gradient-based conditions holds:
\begin{itemize}
\item 
For $l = 0,\ldots,L-1$ the sets $\Bx^*$ and $\Bp^*$ satisfy, 
\begin{equation}\label{def-Forward Propagation}
 \Bx_{l+1}(\Bu^*) = \Lambda_l(\Bx_l(\Bu^*),u^{*}_l),\quad \Bx_0(\Bu^*) \in \X^{B}_0, 
 \end{equation}
\begin{equation}\label{def-backpropagation}
 \Bp_l(\Bu^*) = \partial_{\Bx_l}\Lambda_l\left(\Bx_l(\Bu^*),u^*_l\right)^T\Bp_{l+1}(\Bu^*),\quad 
 \Bp_L(\Bu^*) = - \frac{1}{|B|}\nabla_{\Bx}\Phi_B(\Bx_L(\Bu^*)).
\end{equation}
\item
For $l = 0,\ldots,L-1$, it holds 
\begin{equation}\label{def- Hamiltonian minimization}
\nabla_{u_{l}}J_B(\Bu^*) := \nabla_{u_{l}} \left( \rho R_l(\Bu^*) - \Bp_{l+1}(\Bu^*)\cdot \Lambda_l(\Bx(\Bu^*), u^*_l) \right) = 0,
\end{equation}
\end{itemize}
where the first equality in \eqref{def- Hamiltonian minimization} can be derived by a successive application of the chain rule together with the relation \eqref{def-backpropagation}; See also \cite{Chen2018}. Notice that the adjoint vectors resulting from \eqref{def-backpropagation} are equivalent to the variables derived in the backpropagation process, which motivates the optimal control perspective on this algorithm.

In the framework of optimal control, a similar optimality system in the form of the Pontryagin maximum principle (PMP) can be derived, which is based upon considering  \eqref{Empirical Risk minimization} as a discrete optimal control problem and introducing the following Hamilton-Pontryagin (HP) function
\begin{equation}\label{discrete Hamiltonian}
\mathcal{H}^{B}_l(\Bx_{l}(\Bu),\Bp_{l+1}(\Bu),u_l) := \Bp_{l+1}(\Bu) \cdot \Lambda_l(\Bx_{l}(\Bu),u_l) - \rho R_l(u_l). 
\end{equation}  
With this auxiliary function at hand the system \eqref{def-Forward Propagation}-\eqref{def- Hamiltonian minimization} can equivalently be defined in terms of the HP functions gradients, leading to the observation that for each layer, \eqref{discrete Hamiltonian} is maximized along the optimal trajectory generated by $\Bu^*$. The original statement of the PMP is due to Pontryagin et al. \cite{PBG1962,BGP1956}. One decade after this first formulation of the maximum principle, Halkin \cite{Halkin1966} derived a version of the PMP for discrete-time optimal control problems, which was shortly after further generalized by Holtzmann \cite{Holtzmann1966}. An additional generalization of the PMP was presented by Zhan et al. \cite{Zhan2012} who formulated a unified theory for the maximum principle and weakened the regularity assumptions given by Halkin for the discrete-time case.
Based upon \cite{Zhan2012}, in our framework, the statement of the discrete PMP is as follows:

Suppose that $\Bu^* \in U_{ad}$ solves the supervised learning problem \eqref{Empirical Risk minimization} with corresponding feature maps $\Bx^* := (\Bx_l(\Bu^*))_{l = 0}^{L}$. Then there exist unique backpropagation vectors $\Bp^*:=(\Bp_l(\Bu^*))_{l = 0}^{L}$, such that \eqref{def-Forward Propagation}-\eqref{def-backpropagation} holds and for each layer $l = 0, \ldots, L-1$, the optimal learning parameters satisfy the maximality condition: 
\begin{equation}\label{eHmax}
\cH^{B}_l(\Bx_l(\Bu^*),\Bp_{l+1}(\Bu^*),u^*_l) \geq \cH^{B}_l(\Bx_l(\Bu^*),\Bp_{l+1}(\Bu^*),u), \qquad u \in U_l.
\end{equation}
In its discrete-time version and opposed to the original statement of the PMP in continuous-time, a main assumption is the convexity of the sets $R_l(U_l)$ and $\Lambda_l(\Bx,U_l)$ for any fixed $\Bx \in \X_l^B$. We remark that due to the affine structure of the pre-activation mechanism, we have that in our framework, $\Lambda_l(\Bx,U_l)$ is convex for any fixed $\Bx \in \X_l^B$ if the set $U_l$ is convex. This condition also implies the convexity of $R_l(U_l)$ in the standard CNN framework with a continuous regularizer; see e.g. \cite{LiHaoEMSA2018}. It is important to remark that these convexity assumptions are not required for the construction of our PMP-based algorithm. The maximality condition \eqref{eHmax} allows the formulation of backpropagation algorithms in cases where no gradient or subgradient is available. That is, for example, where the set of admissible parameters is discrete \cite{LiHaoEMSA2018,Baldassi2007,Baldassi2017} or the regularizer is only lower semi-continuous; for the latter, see Section \ref{sec-Sparse network training}.

\section{Generalized training based on the PMP}
	\label{sec-Successive approximation schemes} 

Having a characterization of the optimal feature maps and backpropagation variables in the form of \eqref{def-Forward Propagation} - \eqref{def-backpropagation} together with the maximality condition stated in the discrete PMP, the backpropagation algorithm can be motivated as a forward-backwards sweep method. This class of indirect methods aims to derive an optimal solution $\Bu^{*}$ by first processing forward through the structure \eqref{def-Forward Propagation} utilizing some previously chosen $\Bu^{(k)}$, followed by determining the terminal cost functions gradients, to then propagate them backwards using \eqref{def-backpropagation}. After all the variables $\Bx_l(\Bu^{(k)})$ and $\Bp_l(\Bu^{(k)})$ have been determined, the network's parameters are updated based upon the condition \eqref{def- Hamiltonian minimization}. Corresponding to the classical backpropagation approach, this update is commonly determined with one step of a gradient descent-based strategy. The process is then repeated iteratively, using the new parameters $\Bu^{(k+1)}$. While this classical backpropagation approach is widely applied in machine learning, it is generally not applicable in cases where no gradient or subgradient information with respect to the trainable parameters is available. That is, if the set of admissible parameters is discrete or a discontinuous regularizer is used. However, considering a condition similar to \eqref{eHmax} for updating the network's parameters gives rise to a backpropagation approach which includes gradient-based updating steps as a special case and accommodates training even without available gradient information. The basis for this generalization is the so-called method of successive approximations (MSA), which was introduced by Krylov and Chernous'ko in the early 1960s \cite{Krylov1963}; see also \cite{Kelley1961}. In our discrete-time setting, and in the spirit of the discrete PMP, the MSA is updating the parameters of a CNN by a maximization step on the discrete Hamiltonian as follows
\begin{equation}\label{Def- Hamiltonian maximization}
     u_l^{(k+1)}  = \argmax_{v \in U_l} \; \cH^B_l(\Bx_l(\Bu^{(k)}),\Bp_{l+1}(\Bu^{(k)}),v).
\end{equation}
While fast convergence with these MSA schemes was achieved in special cases, the method appeared less robust with respect to parameter changes and initialization for general nonlinear control systems. As a result, in its basic form, the MSA method is inappropriate for direct application in the framework of machine learning due to its lack of stability. However, in the last decades, several robust variants of the original MSA for continuous-time optimal control problems have been proposed \cite{EWeinanEMSA2018,SakawaShindo1980,ChernouskoLyubushin1982,BreitenbachBorzi2018}, including modifications specifically developed for application to the framework of supervised learning with neural networks  \cite{LiHaoEMSA2018,Hofmann2022}. The sequential quadratic Hamiltonian (SQH) approach for machine learning \cite{Hofmann2022}, which is inspired by its continuous-time equivalent \cite{BreitenbachBorzi2018,BreitenbachBorzi2019,BreitenbachBorzi2020}, enforces the convergence of the basic MSA-scheme by adjusting the Hamiltonian maximization step \eqref{Def- Hamiltonian maximization} using an augmented Hamiltonian approach \cite{SakawaShindo1980,ShindoSakawa1985}. That is, for updating the trainable parameters at each iteration, the classical HP function in \eqref{Def- Hamiltonian maximization} is replaced with the following augmented HP function 
\begin{equation}\label{Batch augmented Hamiltonian}
    \cH^B_{l,\eps}(\Bx_l ({\Bu}),\Bp_{l+1}  ({\Bu}), w_l, u_l) := \cH^B_{l}(\Bx_l ({\Bu}),\Bp_{l+1}  ({\Bu}), w_l) - \frac{\epsilon}{2} \| w_l - u_l\|^2,
\end{equation}
where the augmentation strength $\eps > 0$ is adjusted suitably at each iteration, to guarantee a reduction of the loss functional value. The use of an augmented HP function is motivated by the observation that a decrease in the full batch loss is achieved by increasing the value of the HP function while simultaneously keeping the parameter changes at each iteration moderate; see Theorem \ref{Theorem-sufficient decrease}. While the application of the SQH scheme to the framework of fully-connected ResNet-like architectures is comparably straightforward \cite{Hofmann2022}, much work is required for its extension to the framework of CNNs. In fact, extending the SQH method to CNNs presents significant challenges. One major difficulty is the exact maximization of the augmented HP function, which may be computationally infeasible due to large-scale input data (e.g., high-resolution RGB images) and the complexity of modern CNN architectures with millions of trainable parameters and deep hierarchical structures. To face this challenge, we replace the exact maximization of the augmented Hamiltonian with the relaxed requirement that for all $l = 1,\ldots,L-1$, each update $u_l^{(k+1)} \in U_l$ must only satisfy the following:
\begin{equation}\label{def- ineq augmented Hamiltonian}
\cH^B_{l,\epsilon_k}(\Bx_l ({\Bu^{(k)}}),\Bp_{l+1}  ({\Bu^{(k)}}), u_l^{(k+1)}, u^{(k)}_l)\geq \cH^B_{l}(\Bx_l ({\Bu^{(k)}}),\Bp_{l+1}  ({\Bu^{(k)}}), u^{(k)}_l).
\end{equation}
While this makes the updating step more feasible numerically, processing the complete dataset at once can still lead to longer inference times, increased memory requirements and a costly, potentially prohibitive training process. This motivates a further relaxation of the updating step by invoking a stochastic approach. The basis for this approach is splitting the available dataset $B$ into so-called mini-batches of size $M \leq |B|$ and training at each iteration on one of these randomly sampled subsets. That is, at each iteration, before executing \eqref{def-Forward Propagation} - \eqref{def-backpropagation}, a batch of data is randomly sampled from
\begin{equation}\label{sampling without replacement}
B(M) := \{ B_k \subset 2^{B} \; | \; |B_k| = M \},
\end{equation}
where in the following, we denote the corresponding set of indices for a mini-batch $B_k$ by $I_{B_k}$. The forward-backwards procedure then is executed on this mini-batch of data and the trainable parameters are adjusted based on a mini-batch-specific augmented HP function, where the Hamiltonian is replaced with its mini-batch variant as follows:
\begin{equation}\label{discrete Minibatch-Hamiltonian}
\cH^{B_k}_l(\Bx,\Bp,u) :=  \sum_ {s \in  I_B} \left[p^s \cdot \lambda_l(x^s,u)\right]\chi_{I_{B_k}}(s) - \rho R_l(u), 
\end{equation}  
where $\chi_{I_{B_k}}: I_B \rightarrow \{0,1\}$ defines the indicator function. As a basis for this modification, we establish \eqref{discrete Minibatch-Hamiltonian} as an unbiased estimator for the full-batch HP function \eqref{discrete Hamiltonian}. 

\begin{lemma}\label{Lemma- Unbiased Hamiltonian estimator}
Let $B_k \in B(M)$ be drawn uniformly at random and define $\Bq_l(\Bu):=\frac{|B|}{M}\Bp_l(\Bu)$ for all $l \in \{0,\ldots,L\}$. Then for each fixed $\Bu, \Bw \in U_{ad}$ and all $l = 0, \ldots, L-1$ it holds
\begin{equation*}
 \mathbb{E}\left[\cH^{B_k}_{l}(\Bx_{l}(\Bu),\Bq_{l+1}(\Bu), w_l)\right] = \cH^B_{l}(\Bx_{l}(\Bu),\Bp_{l+1}(\Bu), w_l).  
\end{equation*} 
 
\end{lemma}
\begin{proof}
We assume, without loss of generality for this proof, that $\rho = 0$, as the regularizer is independent of each mini-batch. Since $B_k \in B(M)$ is drawn uniformly at random and due to the fundamental bridge between expectation and probability, we have
\begin{equation}\label{Ch 3: Probability of s in Bk}
\mathbb{E}\left[\chi_{I_{B_k}}\right]  = \mathbb{P}(\chi_{I_{B_k}} = 1) = \mathbb{P}(\{s \in I_{B_k}\}) = \frac{\sum_{B_j \in B(M)} \chi_{{B_j}}(s)}{\binom{M}{|B|}} = \frac{M}{|B|}.
\end{equation}
We reformulate the realizations of the mini-batch HP function as follows
\begin{equation*}\label{Ch 3: Reformulated HP function}
    \cH^{B_k}_l(\Bx_{l}(\Bu),\Bq_{l+1}(\Bu),w_l) = \frac{M}{|B|} \sum_{s\in I_B} \left[q^s_{l+1}(\Bu)\cdot\lambda_l(x^s_{l}(\Bu),w_l)\right]\frac{|B|}{M}\chi_{I_{B_k}}(s). 
\end{equation*}   
Thus, using the result \eqref{Ch 3: Probability of s in Bk}, together with the linearity of the expected value, yields
\begin{align*}
    \mathbb{E}\left[\cH^{B_k}_l(\Bx_{l}(\Bu),\Bq_{l+1}(\Bu),w_l)\right]  
    & = \sum_{s\in I_B} \left[p^{s}_{l+1}(\Bu)\cdot\lambda_l(x^s_{l}(\Bu),w_l)\right]\frac{|B|}{M}\mathbb{E}\left[\chi_{I_{B_k}}\right]\\
    & = \cH_{l}(\Bx_{l}(\Bu),\Bp_{l+1}(\Bu), w_l).
\end{align*}
\end{proof}

With the result of Lemma \ref{Lemma- Unbiased Hamiltonian estimator} at hand, we have that the maximality condition \eqref{eHmax} in the discrete PMP holds in expectation for the mini-batch Hamiltonian \eqref{discrete Minibatch-Hamiltonian}. That is, we have for an optimal $\Bu^{*}\in U_{ad}$ and each layer $l = 0,\ldots,L-1$ that
\begin{equation*}
\mathbb{E}\left[\cH^{B_k}_l(\Bx_l(\Bu^*),\Bq_{l+1}(\Bu^*),u^*_l)\right] \geq \mathbb{E}\left[\cH^{B_k}_l(\Bx_l(\Bu^*),\Bq_{l+1}(\Bu^*),u)\right], \qquad u \in U_l.
\end{equation*}
%
%
This result, together with the law of large numbers, motivates replacing the full-batch Hamiltonian in the augmented HP function \eqref{def- ineq augmented Hamiltonian} by its mini-batch counterpart \eqref{discrete Minibatch-Hamiltonian} and thus gives rise to the following batch SQH (bSQH) scheme: 

\begin{algorithm2e}[H]
 Initialize: \quad $\Bu^{(0)} \in U_{ad}$, \quad $k,j = 0$,\quad $k_{max} \in \N$,\quad $M \leq |B|$ \;
 Hyperparameter: \quad $\epsilon_0 >0$, \quad $\mu>1$, \quad $\eta\in(0,\infty)$\;
 \While{$k < k_{max}$}{
 Draw uniformly at random $B_k \in B(M) $\;
    For $l = 0,...,L-1$\quad  solve\\ \quad $\Bx^{(k)}_{l+1} = \Lambda_l(\Bx^{(k)}_{l},u^{(k)}_l),\quad \Bx^{(k)}_{0} \in \X^{B_k}_{0}$\; 
    For $l = L-1,...,0 $\quad solve\\ \quad
    $ \Bp^{(k)}_l = \partial_{\Bx_l}\Lambda_l\left(\Bx^{(k)}_l,u^{(k)}_l\right)^T\Bp^{(k)}_{l+1}, \quad \Bp^{(k)}_L = -\frac{1}{M}  \nabla_{\Bx}\Phi(\Bx^{(k)}_L)$\;
    For $l = 0,...,L-1$\;
    Choose $w_l \in U_l$ such that
    $\cH^{B_k}_{l,\epsilon_k}(\Bx^{(k)}_l,\Bp^{(k)}_{l+1}, w_l, u^{(k)}_l)\geq \cH^{B_k}_{l}(\Bx^{(k)}_l,\Bp^{(k)}_{l+1}, u^{(k)}_l)$\;
    \While{$J_{B_k}\left(\Bw\right)-J_{B_k}\left(\Bu^{(k)}\right) > -   \eta\normiii{\Bw - \Bu^{(k)}}$}
            {
            Set $j = j + 1, \quad \epsilon_k = \mu^j \hat{\epsilon}_k$\;
            For $l = 0,...,L-1$\;
            Choose $w_l \in U_l$ such that
            $\cH^{B_k}_{l,\epsilon_k}(\Bx^{(k)}_l,\Bp^{(k)}_{l+1}, w_l, u^{(k)}_l)\geq \cH^{B_k}_{l}(\Bx^{(k)}_l,\Bp^{(k)}_{l+1}, u^{(k)}_l)$\;
            }
    Determine $\epsilon_{k+1}=\hat{\epsilon}_{k+1} > 0$ (e.g. \eqref{eqn: SQH},\eqref{eqn: MA})\;
    Set $\Bu^{(k+1)} = \Bw$\;
    Set $ k = k+1, \quad j = 0 $\;
      
}
\caption{bSQH}
\label{algoSQH}
\end{algorithm2e}

The bSQH algorithm adjusts at each iteration the initial guess $\hat{\epsilon}_k > 0$ corresponding to
\begin{equation}\label{Successive increase aug. param.}
\epsilon_{k} = \mu^{j} \hat{\epsilon}_{k},
\end{equation}
by choosing the minimal $j \in \N$ such that a decrease in the mini-batch loss is achieved. It is apparent that the choice of $\hat{\eps}_{k}$ has a decisive impact on the number of required inference steps and thus on the efficiency of the scheme. While the basic SQH scheme uses a backtracking approach as follows 
\begin{equation}\label{eqn: SQH}
\hat{\eps}_{k+1} = \zeta \eps_k \tag{\text{SQH$_{\zeta}$}},
\end{equation}
we propose choosing the initial guess for the augmentation parameter based on the information encoded in the realizations $\{\epsilon_i\}_{i = 1,\ldots,k-1}$ to reduce computational overhead. This gives rise to the following modifications inspired by smoothing techniques known from time-series analysis:
\begin{equation}\label{eqn: MA}
\hat{\epsilon}_{k+1} = 
\begin{cases}
    \zeta \hat{\epsilon}_{k} & \text{if } j = 0 \\
    \frac{1}{\min(k+1, \omega+1)} \sum_{i = k - \min(k, \omega)}^{k} \epsilon_i & \text{else}
\end{cases},
\tag{\text{bSQH$_{\omega}$}} 
\end{equation}
where $\zeta \in (0,1]$ and $\omega \in \N$, must be chosen a-priori. The basis for this modification is the assumption that the choice of $\eps_k$ at each iteration implicitly depends on the local Lipschitz constants of the current mini-batch HP function. Thus, by taking the moving average of the approximations $\eps_k$ for these unknown constants, one expects to achieve an estimate for an optimal augmentation parameter in the next iteration. Even tough we do not give a rigorous proof for the effectiveness of our adjustment \eqref{eqn: MA}, we provide numerical evidence in Section \ref{sec-numerical experiments} which shows that this adjustment can significantly reduce the number of forward steps required to determine a suitable augmentation parameter $\eps_k$ at each iteration. Moreover, while a suitable choice of the augmentation parameter can reduce the number of forward steps, introducing sparsity into the parameters of a CNN can reduce both the inference times and the memory storage of a trained model.     

\section{Sparse network training based on the bSQH scheme}
	\label{sec-Sparse network training} 

The most direct path to introduce sparsity into the solutions to  \eqref{Empirical Risk minimization} is a penalty term in the form of the $L^0$-'norm' which, for a given vector $v \in \R^n$, is defined by
\begin{equation}\label{L0 semi-norm}
\|v\|_0 :=|\{i \in \{1,\ldots ,n\}\;|\; v_i \neq 0 \}|,
\end{equation}
and thus strongly penalizes the number of non-zero parameters in a CNN. 
However, while this regularization technique is the most direct to introduce sparsity, it is seldom used in machine learning, due to its non-differentiability and optimization complexity, which makes classical backpropagation infeasible. The formulation of the updating step in the bSQH approach on the contrary, allows a straightforward derivation of an explicit update even in this non-continuous case. To prove this conjecture, we consider non-residual network structures similar to the well-known LeNet-5 architecture, where a linear output layer is used. We remark that a linear output layer can always be assumed without loss of generality by incorporating the non-linearity of the last layer into the terminal loss function. With the notation of Proposition \ref{Proposition-boundedness of x}, the corresponding forward propagation of such a network can be described as follows
\begin{equation}\label{MLP with linear layer}
\begin{aligned}
	x_{L} & = \cW_{L-1}(u_{L-1}) x_{L-1} + \cB_{L-1}(u_{L-1}),\\
    x_{l+1} &= f_l(\cW_l(u_{l}) x_{l} + \cB_l(u_{l})) , \quad l = 0,\ldots, L-2,\\
    x_0 &= x^s.
\end{aligned}
\end{equation}
The composite structure of \eqref{MLP with linear layer} enables the reordering of the nonlinearity $f_l$, and gives rise to a network with layer functions defined by

\begin{equation}\label{Control affine MLP}
\lambda_l(x,u) = 
\begin{cases}
      \cW_l(u)f_{l-1}(x) + \cB_lu, \quad \text{for} \quad l = 1, \ldots,L-1\\
\cW_{l}(u) x + \cB_l u, \quad \qquad \text{for} \quad l = 0.  
\end{cases}
\end{equation}
An important observation for our approach is that for a given set of parameters $\Bu$, both the network defined by \eqref{MLP with linear layer} and the one defined corresponding to \eqref{Control affine MLP} have the same output $x_L$. Thus, for architectures such as \eqref{MLP with linear layer}, we can formulate the supervised learning task \eqref{Empirical Risk minimization} equivalently with neural networks constraints that admit layer architectures given by \eqref{Control affine MLP}. Due to the control-affine structure of \eqref{Control affine MLP}, this in turn gives rise to the following augmented Hamiltonian function
$$
\cH^{B_k}_{l,\eps}(\Bx ,\Bp , w, u) = \cF^{B_k}_l(\Bx,\Bp)u - \left( \rho R_l (u)+ \frac{\epsilon}{2} \| w - u\|^2 \right),
\quad \cF^{B_k}_l(\Bx,\Bp) : = \sum_{s \in I_{B_k}}p^{s}\cdot(F_l(x^{s}) + \cB_l),
$$

where $F_l(x^s)$, for a fixed state $x^s$, denotes a suitable transformation matrix of the linear transformations $u \mapsto \cW_l(u)f_{l-1}(x^s)$ and $u \mapsto \cW_0(u)x^s$ respectively. This formulation of the augmented Hamiltonian at hand, the updating step in the bSQH algorithm, simplifies to the following regularized quadratic problem
\begin{equation}\label{proximal subproblem}
\begin{aligned}
  u_l^{(k+1)} &= \argmax_{u \in U_l} \; \mathcal{H}^{B_k}_{l,\eps_k}(\Bx_l^{(k)},\Bp_{l+1}^{(k)}, u, u^{(k)}_l)\\
  & = \argmin_{u \in U_l} \frac{1}{2} \left\| u - \left(w^{(k)}_l +  \frac{1}{\epsilon_k} \cF^{B_k}_l(\Bx^{(k)}_l,\Bp^{(k)}_{l+1}) \right) \right\|^2 + \frac{\rho}{\epsilon_k} R_l(u).
\end{aligned}
\end{equation}
With a suitable regularizing function $R_l$, the optimization problem \eqref{proximal subproblem} admits a unique and analytic solution. Based on this reformulation and to enforce sparsity of a trained network by simultaneously increasing the stability of the forward propagation \cite{Haber2017}, we propose a regularization approach similar to the well-known Elastic-Net regularization \cite{Zou2005}. Specifically, we replace the $L^1$-norm in the  Elastic-Net approach with \eqref{L0 semi-norm} to obtain the following regularizer 
\begin{equation}\label{L2L0 regularization}
    R_l(u) = \frac{\alpha}{2} \|u\|^2 + (1 - \alpha)\|u\|_0, \qquad 
\end{equation}
for all $l = 0,\ldots,L-1$, where $\alpha \in [0,1)$. 
Due to the formulation \eqref{proximal subproblem}, each update of the trainable parameters at each iteration of Algorithm \ref{algoSQH} takes for all $l = 0,\ldots, L-1$ the following explicit form
\begin{equation}\label{proximal subproblem L0L2}
\begin{aligned}
  u_l^{(k+1)} &= \argmin_{\Bu \in U_l} \frac{1}{2} \left\|u - \frac{\epsilon_k}{\epsilon_k + \alpha\rho}\left(u^{(k)}_l + \frac{1}{\epsilon_k} \cF^{B_k}_l(\Bx^{(k)}_l,\Bp^{(k)}_{l+1})\right) \right\|^2 + \frac{(1 -\alpha)\rho}{\epsilon_k + \alpha\rho}\|u\|_0\\
  &= \cT_{\frac{(1 -\alpha)\rho}{\epsilon_k + \alpha\rho}}\left(\frac{\epsilon_k}{\epsilon_k + \alpha\rho}\left(u^{(k)}_l + \frac{1}{\epsilon_k} \cF^{B_k}_l(\Bx^{(k)}_l,\Bp^{(k)}_{l+1})\right)\right),
\end{aligned}
\end{equation}
where the Hard-threshold operator for a given vector $v \in \R^n$ and $\gamma > 0$ is defined by 
$$
\cT_{\gamma} (v)_i=
\begin{cases} 
0 & |v_i| \leq \sqrt{2\gamma} \\ 
v_i & \text{otherwise}\end{cases}, \qquad i = 1, \ldots, n.
$$

Notice that a similar reasoning applies in the Elastic-Net approach, where the $\cT$ in \eqref{proximal subproblem L0L2} is replaced with the Soft-threshold operator defined as follows for a given vector $v \in \R^n$ and $\gamma > 0$:
$$
\cS_{\gamma} (v)_i=
\begin{cases} 
0 & |v_i| \leq \gamma \\ 
\left(v_i - \gamma \right)\sgn(v_i) & \text{otherwise}\end{cases}, \qquad i = 1, \ldots, n.
$$
However, while the latter penalty only encourages sparsity of the trained network's parameters, a penalty of the form \eqref{L2L0 regularization} enforces exact sparsity of a trained CNN. We present a comparison of the sparsity-inducing properties of both approaches in Section \ref{sec-numerical experiments}.

	\section{Convergence of the bSQH scheme}
 \label{Sec-Well-definedness and Convergence}
 \vspace{0.3cm}

In the following, we present results which prove that each updating step of the bSQH algorithm is feasible and can be generated within a finite number of steps.  We remark that while some of these upcoming results are similar to those stated in \cite{Hofmann2022}, we present here a rigorous proof for the validity of the sufficient decrease condition
\begin{equation}\label{Ch 4.2: Sufficient decrease condition}
J_{B_k}\left(\Bu^{(k+1)}\right)-J_{B_k}\left(\Bu^{(k)}\right) \leq -   \eta \normiii{ \Bu^{(k+1)} - \Bu^{(k)}},
\end{equation}
under the relaxed assumptions \ref{A1} - \ref{A2}, which include those discussed in \cite{Hofmann2022} as a special case. As a basis for the derivation of these results, we give the following auxiliary lemma.

\begin{lemma}\label{Lemma- boundedness of p}
    Let \ref{A1} - \ref{A2} hold and define $\Bu \mapsto \Bx_l(\Bu)$ and $\Bu \mapsto\Bp_l(\Bu)$ for each $\Bu\in U_{ad}$ recursively by \eqref{def-Forward Propagation} and \eqref{def-backpropagation} respectively. Then for all $l = 0,\ldots,L-1$. the control to adjoint map $\Bp_l(\;\cdot\;)$ is uniformly bounded and there exists a constant $C>0$ such that for any $\Bu,\Bw \in U_{ad}$ it holds
    \begin{equation}\label{Lipschitz property control-to-feature map}
        \|\Bx_l(\Bu) - \Bx_l(\Bw) \| \leq C \sum_{j = 0}^{L-1} \| u_j - w_j \|,
        \qquad  
        \|\Bp_l(\Bu) - \Bp_l(\Bw) \| \leq C \sum_{j = 0}^{L-1} \| u_j - w_j \|,
    \end{equation}
    for all $l = 0,\ldots, L-1$.
\end{lemma}
\begin{proof}
Let $\Bu,\Bw \in U_{ad}$, then due to the Lipschitz continuity of $\Lambda_l$ for all $l \in \{0,\ldots,L-1\}$, by assumption \ref{A2}, there exists a constant $C>0$ such that
\begin{equation}
\| \Xp(\Bu)  -\Xp(\Bw) \| \leq C \left( \| \Xl(\Bu)  -\Xl(\Bw)\| +  \| u_l - w_l \| \right).
\end{equation}
Then, by a discrete version of the Gronwall lemma, we have a constant $C > 0$ such that for all $l = 0,\ldots,L-1$, it holds
\begin{equation*}
\| \Bx_l(\Bu)  -\Bx_l(\Bw) \| \leq  C \sum_{j = 0}^{l-1} \| u_j - w_j \|, 
\end{equation*}
which proves the Lipschitz property \eqref{Lipschitz property control-to-feature map} of the control-to-feature maps. For the adjoint map, a similar reasoning applies.
To prove the boundedness of the backpropagation variables, we have by the boundedness of $\partial_x \Lambda_l$, that there exists some constant $C > 0$ for all $l = 0,\ldots,L-1$ and any $\Bu \in U_{ad}$, such that
\begin{equation*}
\| \Pl(\Bu) \| =  \| \partial_x \Lambda_l(\Xl(\Bu),u_l)^T  \Pp(\Bu)  \| \leq C\| \Pp(\Bu) \|.
\end{equation*}
Thus, by a change of index and an application of the discrete Gronwall lemma, there exists a constant $C>0$ such that the following holds for all $\Bu \in U_{ad}$
\begin{equation*}
\| \Pl(\Bu) \| \leq  C^L\| \Bp_L(\Bu) \| = C^L\| \nabla \Phi_{B}(\Bx_L(\Bu)) \|.
\end{equation*}
The boundedness of $\nabla \Phi_{B}$ by assumption \ref{A1} now proves the uniform boundedness of the parameter to adjoint maps.
\end{proof}
Next, we prove that a successful update with the bSQH algorithm can be achieved within a finite number of steps. That is, given the iterate $\Bu^{(k)}$ and $\eta\in (0,\infty)$, at each iteration a suitable set of parameters $\epsilon_k>0$ exists such that the sufficient decrease condition \eqref{Ch 4.2: Sufficient decrease condition} holds.
\begin{theorem}\label{Theorem-sufficient decrease}
Let $B_k \in B(M)$ and assume that \ref{A1}-\ref{A2} hold. Then, there exits a constant $C > 0$ independent of $\epsilon_k > 0$, such that for every $\Bu, \Bw \in U_{ad}$ satisfying
\begin{equation}\label{Ch4.2.1 augmented Hamiltonian inequality}
\cH^{B_k}_{l,\epsilon_k}(\Bx_l ({\Bu}),\Bp_{l+1} ({\Bu}), w_l, u_l)\geq \cH^{B_k}_{l}(\Bx_l ({\Bu}),\Bp_{l+1}  ({\Bu}), u_l), \qquad l = 1,\ldots,L-1,
\end{equation}
it holds
\begin{equation}\label{Ch4.2.1 Sufficient decrease condition}
J_{B_k}(\Bw)-J_{B_k}(\Bu) \leq (C - \epsilon_k) \normiii{ \Bw - \Bu }.
\end{equation}
\end{theorem}
\begin{proof}
For ease of notation, assume $|B_k| = 1$ and suppress the parameters $s$ and $B_k$. It is worth noting that this adjustment can be made without loss of generality, due to the summable structure of the HP and cost function. By definition, for any $\Bw, \Bu \in U_{ad}$, we have:
\begin{equation}\label{Ch.4.2.1: Theorem Estimate 1}
\begin{aligned}
J(\Bw) - J(\Bu) &= \Phi(\Bx_{L}(\Bw))  + \sum_{l = 0}^{L-1} R_l( w_l )-  \left( \Phi(\Bx_{L}(\Bu)) + \sum_{l = 0}^{L-1} R_l( u_l )\right) \\
& =  \Phi(\Bx_{L}(\Bw))  - \Phi(\Bx_{L}(\Bu))  - \Delta \cH (\Bw,\Bu) \\
&\quad + \sum_{l = 0}^{L-1} \Pp(\Bu) \cdot \left( \Xp(\Bw) - \Xp(\Bu) \right)\\
&\quad - \sum_{l = 0}^{L-1} \cH_l(\Xl(\Bw),\Pp(\Bu), w_l) - \cH_l(\Xl(\Bu),\Pp(\Bu), w_l),
\end{aligned}
\end{equation}  
where $\Delta \cH (\Bw,\Bu) : =  \sum_{l = 0}^{L-1} \cH_l(\Xl({\Bu}),\Pp(\Bu),w_l) - \cH_l(\Xl({\Bu}),\Pp(\Bu),u_l)$. For the following, we define for all $l = 0,\ldots,L-1$, the convex combination $x_l(\tau) = \Xl(\Bw) + \tau(\Xl(\Bu) - \Xl(\Bw))$, where $\tau \in (0,1)$. Due to the convexity of $\X^B_l$ for each $l = 0,\ldots, L-1 $, the mean-value theorem guarantees the existence of  $ \tau_l \in (0,1) $ such that 
\begin{equation}\label{Ch.4.2.1: Theorem Equality 1}
\begin{aligned}
&\cH_l(\Xl(\Bw),\Pp(\Bu), w_l) - \cH_l(\Xl(\Bu),\Pp(\Bu), w_l)\\
& = \left\langle \nabla_x \cH_l(x_l(\tau_l),\Pp(\Bu), w_l), \Xl(\Bu) - \Xl(\Bw) \right\rangle \\
& = \langle \Pp(\Bu)^T \left(\partial_x \Lambda_l(x_l(\tau_l),w_l) - \partial_x \Lambda_l(\Xl(\Bu),u_l)\right), \left(\Xl(\Bu) - \Xl(\Bw)\right)\rangle \\
&\qquad + \langle \Pl(\Bu), \left(\Xl(\Bu) - \Xl(\Bw)\right)\rangle.
\end{aligned}
\end{equation}
Similarly, we have that there exists $ \tilde{\tau} \in (0,1) $ such that 
\begin{equation}\label{Ch.4.2.1: Theorem Equality 2}
\Phi(\Bx_{L}(\Bw)) - \Phi(\Bx_{L}(\Bu)) = \langle \nabla\Phi(x_{L}(\tilde{\tau})),\Bx_L(\Bu) - \Bx_L(\Bw)\rangle.
\end{equation}
Substituting both \eqref{Ch.4.2.1: Theorem Equality 1} and \eqref{Ch.4.2.1: Theorem Equality 2} into \eqref{Ch.4.2.1: Theorem Estimate 1}, resolving the involved telescopic sum and using $\Bp_L(\Bu) = \nabla\Phi(\Bx_{L}(\Bu)) $, we achieve 
\begin{equation}\label{Ch.4.2.1: Theorem Estimate 2}
\begin{aligned}
J(\Bw) - J(\Bu) &=   - \Delta \cH (\Bw,\Bu)\\
&\quad + \langle \nabla\Phi(x_{L}(\tilde{\tau})) - \nabla\Phi(\Bx_{L}(\Bu)), \Bx_L(\Bw) - \Bx_L(\Bu)\rangle \\
&\quad + \sum_{l=0}^{L-1} \langle( \Pp(\Bu)^T\left(\partial_x \Lambda_l(x_l(\tau_l),w_l) - \partial_x \Lambda_l(\Xl(\Bu),u_l)\right),  \Xl(\Bw) - \Xl(\Bu)\rangle.
\end{aligned}
\end{equation}  
Using the Cauchy-Schwarz inequality and the Lipschitz continuity of $\nabla \Phi$, there exists a constant $C>0$ such that 
\begin{equation}\label{Ch.4.2.1: Theorem Estimate 3}
 \langle \nabla\Phi(x_{L}(\tilde{\tau})) - \nabla\Phi(\Bx_{L}(\Bu)),\Bx_L(\Bw) - \Bx_L(\Bu)\rangle  \leq  (C+\tilde{\tau})\| \Bx_{L}(\Bw)  - \Bx_{L}(\Bu))\|^2. 
\end{equation}
Similarly by the Lipschitz continuity of each $\partial_x \Lambda_l$ uniformly on $\X^B_l \times U_l$, there exists a constant $C>0$ such that 
\begin{equation}\label{Ch.4.2.1: Theorem Estimate 4}
\begin{aligned}
& \langle( \Pp(\Bu)^T\left(\partial_x \Lambda_l(x_l(\tau_l),w_l) - \partial_x \Lambda_l(\Xl(\Bu),u_l)\right),  \Xl(\Bw) - \Xl(\Bu)\rangle\\
& \qquad \leq  \| \Pp(\Bu)\| \| \partial_x \Lambda_l(x_l(\tau_l),w_l) - \partial_x \Lambda_l(\Xl(\Bu),u_l)\|\| \left(\Xl(\Bw) - \Xl(\Bu)\right)\|\\
& \qquad \leq  \| \Pp(\Bu)\| \left(C \| w_l -u_l \| \| \Xl(\Bw) - \Xl(\Bu) \| +(C+\tau_l)\| \Xl(\Bw) - \Xl(\Bu) \|^2\right)
\end{aligned}
\end{equation}
Substituting the estimates of \eqref{Ch.4.2.1: Theorem Estimate 3} and \eqref{Ch.4.2.1: Theorem Estimate 4} into \eqref{Ch.4.2.1: Theorem Estimate 2} and using Lemma \ref{Lemma- boundedness of p} together with the Cauchy-Schwarz inequality, guarantees the existence of a constant $C>0$ such that  for all $\Bw,\Bu \in U_{ad}$ we have
\begin{equation}\label{Ch 4.2.1: Cost decrease}
J(\Bw) - J(\Bu) \leq   - \Delta \cH (\Bw,\Bu) +  C \normiii{\Bw  -\Bu }.   
\end{equation}
Suppose that $\Bw, \Bu  \in U_{ad}$ are chosen such that  \eqref{Ch4.2.1 augmented Hamiltonian inequality} holds. Then, by summing up over all layers, it holds: 
\begin{equation}
\Delta \cH (\Bw,\Bu) \geq \epsilon_k \normiii{\Bw - \Bu}.
\end{equation}
By substituting this into \eqref{Ch 4.2.1: Cost decrease} we achieve the desired result, where the constant $C$ is independent of $\epsilon_k$.
\end{proof}
As a result of Theorem \ref{Theorem-sufficient decrease}, the sufficient decrease condition \eqref{Ch 4.2: Sufficient decrease condition} on the batch individual loss is guaranteed as soon as an augmentation parameter is determined satisfying $ \epsilon_k > \eta + C$. Notice that this is ensured by the successive increase of the augmentation parameter $\epsilon_k$ corresponding to \eqref{Successive increase aug. param.}. In the particular case of the full-batch setting of Algorithm \ref{algoSQH}, this leads to a monotonically decreasing sequence of cost functional values.
\begin{lemma}\label{Lemma-Convergence}
Suppose $B_k = B$  for all $k \in \N$ and let the sequence  $(\Bu^{(k)})_{k \in \N} \subset U_{ad}$ be defined such that $J_B(\Bu^{(0)}) < \infty$ and the sufficient decrease condition  \eqref{Ch 4.2: Sufficient decrease condition} holds for all $k \in \N$. Then if $J_B$ is bounded from below, the sequence $(J_B(\Bu^{(k)}))_{k\in \mathbb{N}}$ of corresponding cost function values is monotonically decreasing and for $(\Bu^{(k)})_{k \in \N}$ it holds $
\lim_{k \rightarrow \infty } \normiii{ \Bu^{(k)} - \Bu^{(k+1)} } = 0. 
$.
\end{lemma}

\begin{proof}
As the sufficient decrease condition \eqref{Ch 4.2: Sufficient decrease condition} holds for each $k \in \N$, the sequence $(J_B(\Bu^{(k)}))_{k\in \mathbb{N}}$ is monotonically decreasing and the boundedness from below of $J_B$ implies convergence.  Summing  \eqref{Ch 4.2: Sufficient decrease condition} over the first $m \in \N$ members of the corresponding sequences, we obtain
$
J_{B}(\Bu^{(0)}) -   J_{B}(\Bu^{(m)}) >  \sum_{k =0}^{m}  \eta\normiii{\Bu^{(k+1)}- \Bu^{(k)} }.
$
Due to $J_B(\Bu^{(0)}) < \infty$ and by the convergence of the sequence of cost function values, we deduce that $m \rightarrow \infty$ implies 
$
\lim_{k \rightarrow \infty } \normiii{ \Bu^{(k)} - \Bu^{(k+1)} } = 0. 
$
\end{proof}

We remark that for both of the prior results, only the Assumptions \ref{A1}-\ref{A2} and the boundedness from below of $J_B$ were required. Therefore, the bSQH scheme is also appropriate for an application to supervised learning problems, where the regularizer is discontinuous and given by \eqref{L2L0 regularization}. Further, opposed to the discrete PMP formulation, the convexity of sets $U_l$ is not a vital requirement for the applicability of Algorithm \ref{algoSQH} either. As a result, and in contrast to the classical backpropagation framework, our methodology could readily be used to train discrete-weight neural networks. 

However, we can further extend these results if the HP function admits a generalized directional derivative with respect to the trainable parameters and the maximization of the augmented HP function holds at least asymptotically.
\begin{theorem}\label{Theorem-Convergence bSQH}
Let \ref{A1}-\ref{A2} hold and suppose that for all $l = 0, \ldots, L-1$ the functions $(x,u) \mapsto \partial_u \Lambda_{l}(x,u)$ exist and are jointly-continuous.
Further, let $(\Bu^{(k)})_{k \in \N} \subset U_{ad}$ define a bounded sequence with corresponding feature maps and back-propagation variables $\By^{(k)}_l := (\Xl(\Bu^{(k)}),\Pp(\Bu^{(k)}))$ generated by Algorithm \ref{algoSQH}. If the corresponding regularizers $R_l: U_l \rightarrow \R$ are locally Lipschitz and there exists a non-zero, null sequence $(e^{(k)})_{k \in \N}$ such that 
\begin{equation}\label{Ch 4.2: almost optimal iterate}
\sum_{l = 0}^{L-1} \max_{w \in U_l} \; \cH^B_{l,\eps_k}(\By_l^{(k)}, w, u^{(k)}_l)  - \cH^B_{l,\eps_k}(\By_l^{(k)},u^{(k+1)}_l, u^{(k)}_l)  \leq e^{(k)},
\end{equation}
then every accumulation point $\bar{\Bu} \in U_{ad}$ of $(\Bu^{(k)})_{k \in \N}$ satisfies  
    \begin{equation}\label{Ch4.2: Optimality of Acc Point}
        0 \in \partial_{u_l}^CJ_B(\bar{\Bu}),
    \end{equation}
for all $l = 0, \ldots, L-1$.
\end{theorem}
\begin{proof}
Since all remaining cases can be derived simultaneously, we suppose $|B| = 1$ and suppress the occurrence of each $s$ and $B$ for ease of notation. Further, for all $l = 0,\ldots,L-1$ we define
    $$
    h^{(k)}_{l}(v) := - \cH_{l,\eps_k}(\By_l^{(k)}, v,u^{(k)}_l).
    $$
Due to the boundedness of $(\Bu^{(k)})_{k \in \N} \subset U_{ad}$, there exists a convergent subsequence (also denoted $(\Bu^{(k)})_{k \in \N}$ ) with an unique accumulation point $\bar{\Bu} \in U_{ad}$. Due to \eqref{Ch 4.2: almost optimal iterate} the variational principle of Ivar Ekeland \cite{Ekeland1974} guarantees that for every $u_l^{(k+1)} \in U_l $ there exists some $w_l^{(k+1)} \in U_l $ such that 
\begin{equation}\label{Ch 4.2: Ekeland convergence}
       \|w_l^{(k+1)} - u_l^{(k+1)}\| \leq \sqrt{e^{(k)}},
\end{equation}
and
\begin{equation}\label{Ch 4.2: Ekeland inequality}
        h^{(k)}_{l}(u_l) \geq h^{(k)}_{l}(w_l^{(k+1)}) - \sqrt{e^{(k)}}\|w_l^{(k+1)} - u_l\|,\qquad \forall u_l \in U_l. 
\end{equation}
By setting $u_l = w_l^{(k+1)} + t v_l$ for $t \geq 0$, we achieve equivalently to \eqref{Ch 4.2: Ekeland inequality} for all $k \in \N$
\begin{equation}\label{Ch 4.2: Ekeland estimate}
    \frac{h^{(k)}_{l}(w_l^{(k+1)} + t v_l) - h^{(k)}_{l}(w_l^{(k+1)}) }{t}  \geq -\sqrt{e^{(k)}}\|v_l\|, \qquad \forall v_l \in U_l,\; t \geq 0.
\end{equation}
We apply the definition of the generalized directional derivative \eqref{generalized directional derivative} to \eqref{Ch 4.2: Ekeland estimate} and use the sum rule of the generalized directional derivative \cite{Clarke2013} together with the Cauchy-Schwarz inequality to derive for all $v_l \in U_l$:
\begin{equation}\label{Ch 4.2: directional derivative lower bound}
   \langle \partial_u\Lambda_l(\Bx_l(\Bu^{(k)}),u_l^{(k)}))^T\Bp_{l+1}(\Bu^{(k)}), v_l \rangle + \mathring{R}_l(w_l^{(k+1)}; v_l) \geq - \left( \sqrt{e^{(k)}} + \eps_{k}\|w_l^{(k+1)} -u_l^{(k)} \| \right) \|v_l\|.
\end{equation}
Notice that by Theorem \ref{Theorem-sufficient decrease}, we can guarantee that there exists a constant $C>0$ such that for all augmentation parameters $\eps_k>0$ generated in the process of Algorithm \ref{algoSQH}, it holds $\eps_k \leq \Bar{\eps}: = \mu(\eta + C) $. Then, by taking the limes superior with respect to $k$, we achieve for the right-hand side of \eqref{Ch 4.2: directional derivative lower bound}, due to \eqref{Ch 4.2: Ekeland convergence}, as well as $e^{(k)} \rightarrow 0$ and the convergence of $u_l^{(k)}$ as $k \rightarrow \infty$ the following
\begin{equation*}
\begin{aligned}
0 &\leq \limsup_{k \rightarrow \infty} \left( \sqrt{e^{(k)}} + \eps_k\|w_l^{(k+1)} -u_l^{(k)} \| \right)\\
  & \leq  \lim_{k \rightarrow \infty} (1 + \bar{\epsilon})\sqrt{e^{(k)}} +  \lim_{k \rightarrow \infty} \bar{\epsilon} \|u_l^{(k+1)} -u_l^{(k)} \| = 0.
\end{aligned}
\end{equation*}
For deriving the limit for the left-hand side of \eqref{Ch 4.2: directional derivative lower bound}, we take into account the upper semicontinuity of the generalized directional derivative \cite{Clarke2013}, the continuity of $\partial_u\lambda_l$ and the Lipschitz properties of $\Bu \mapsto \Bx_l(\Bu)$ and $\Bu \mapsto \Bp_l(\Bu)$ for all $l = 0,\ldots,L-1$ by Lemma \ref{Lemma- boundedness of p}, to derive
\begin{equation}\label{Ch 4.2: Limsup Hamiltonian}
\begin{aligned}
   & \limsup_{k \rightarrow \infty} \left( \left\langle \partial_u\Lambda_l(\Bx_l(\Bu^{(k)}),u_l^{(k)}))^T\Bp_{l+1}(\Bu^{(k)}), v_l \right\rangle + \mathring{R}_l(w_l^{(k+1)}; v_l)\right)\\
   &\qquad \qquad \leq \left( \left\langle \partial_u\Lambda_l(\Bx_l(\bar{\Bu}),\bar{u}_l))^T\Bp_{l+1}(\bar{\Bu}), v_l \right\rangle + \mathring{R}_l(\bar{u}_l; v_l)\right).     
\end{aligned}
\end{equation}
This, together with the sum rule for the Clarke subdifferential \cite{Clarke2013}, implies 
$$
0 \in \partial^C \left[\cH^{B}_l(\Bp_{l+1}(\bar{\Bu}),\Bx_{l}(\bar{\Bu}), \;\cdot\;)\right](\bar{u}_l),
$$
for all $l =0,\ldots,L-1$. The equivalence to \eqref{Ch4.2: Optimality of Acc Point} is then established inductively over $l =0,\ldots,L-1$ by applying the relation \eqref{def-backpropagation} and starting from
$$
\Bp_{L-1}(\bar{\Bu}) = \partial_{u_{L-1}}\Lambda_l(\Bx_{L-1}(\bar{\Bu}),\bar{u}_{L-1}))^T\Bp_{L}(\bar{\Bu}) = \nabla_{u_{L-1}} \Phi(\Bx_{L}(\bar{\Bu})),
$$
which concludes the proof.
\end{proof}
The conditions of Theorem \ref{Theorem-Convergence bSQH} are easily satisfied in the full-batch setting due to Lemma \ref{Lemma-Convergence} and when using a reformulated network such as \eqref{Control affine MLP} together with an updating strategy corresponding to \eqref{proximal subproblem}. While this theorem applies to supervised learning problems with Elastic-Net regularization, it also extends to other regularization approaches. Without considering a reformulated network, the maximization of the Hamiltonian function becomes more complex as the activation occurs after the affine transformation. Nevertheless, the property \eqref{Ch 4.2: almost optimal iterate} can be enforced by taking into account that for $\epsilon_k$ sufficiently large, the augmented Hamiltonian becomes strictly concave \cite{Stachurski1989}. Thus, methods such as quasi-Newton methods can determine a close approximation of the true global maximum of the augmented Hamiltonian within a few iterations and thereby enable convergence of the full-batch bSQH in these cases as well. In the mini-batch framework and due to the results of Lemma \ref{Lemma- Unbiased Hamiltonian estimator}, the 
property \eqref{Ch 4.2: almost optimal iterate} can be enforced at least asymptotically by successively increasing the mini-batch size throughout training. We will validate numerically that the property \eqref{Ch 4.2: almost optimal iterate} can similarly be enforced for a fixed mini-batch size by combining it with the approach \eqref{eqn: MA}, where $\zeta =1$. 

	\section{Numerical experiments}
	\label{sec-numerical experiments}
\setkeys{Gin}{draft = False}
In the following, we numerically investigate the bSQH algorithm as a method for training sparse CNNs. A  measure of sparsity for a set of parameters $\Bu \in U_{ad}$ is provided by the ratio
$$
\text{Sp}(\%) :=  100\left( 1 - \frac{\sum_{l=0}^{L-1} \|u_l\|_{0}}{\sum_{l=0}^{L-1}{m_l}} \right),
$$
where $m_l$ is the dimension of the underlying parameter space $U_l$.
This investigation is split into three parts. First, the full-batch variant of Algorithm \ref{algoSQH} is examined with both the sparsity-enforcing Elastic-Net regularization and the non-continuous $L^0$-'norm' based regularizier \eqref{L2L0 regularization}. In both cases, we also investigate the impact of the selection strategies \eqref{eqn: SQH} and \eqref{eqn: MA} on the computational overhead of the bSQH scheme. The second part of our numerical investigation is concerned with the mini-batch variant of Algorithm \ref{algoSQH}. In particular, we investigate the impact of the mini-batch approach on the stability of the algorithm's convergence and demonstrate the computational advantages of the mini-batch approach compared to its full-batch counterpart. In the third and last part of this section, we present an application of the bSQH method in its mini-batch setting to a multi-classification task for medical imaging. These experiments substantiate our claim that the results of this work still hold in the case of network architectures using ReLU activation functions and max-pooling operations, where the differentiability assumptions in \ref{A2} do not apply. In all of the upcoming experiments, the weights of the network architectures are initialized using Xavier-normal initialization \cite{Xavier2010}, while at the start of training, the bias vectors are set to zero. An appropriate parameter update $w\in U_l$ for every layer at each iteration of the bSQH scheme is generated corresponding to the updating strategy \eqref{proximal subproblem} which readily applies in the framework of a LeNet-5 architecture and guarantees that the property \eqref{Ch 4.2: almost optimal iterate} is satisfied exactly in the full-batch framework. This reformulation-based updating approach is different to the L-BFGS-based \cite{Liu1989} updating technique commonly used in the PMP-based machine learning framework; e.g. see \cite{EWeinanEMSA2018,Hofmann2022} and enables training neural networks with non-continuous regularizers such as the $L^0$-'norm', which directly penalizes the number of non-zero parameters of a network. 

\subsection{Sparse network training}
For our investigation in the full-batch setting, we train a standard LeNet-5 architecture to classify handwritten digits in the MNIST dataset \cite{deng2012mnist}. The network structure is chosen as proposed in \cite{LeCun1988} and uses $\tanh$ activations together with average pooling layers. In our experiments, we consider cross-entropy loss and test our algorithm for both of the sparsity-enforcing regularizers, Elastic-Net and \eqref{L2L0 regularization}, where  $\alpha = 0.8$ and $\rho = 1e-04$.

The MNIST dataset consists of $60.000$ images for testing and $10.000$ images for training. For our numerical investigation, we extract a suitable subset of $20.000$ images from the testing dataset and apply data normalization as a preprocessing step. The LeNet-5 is then trained in a full-batch fashion using a maximum of $k_{max} = 1500$ iterations of Algorithm \ref{algoSQH}, where we set $\eta=10e-10$ and initialize the augmentation parameter with $\epsilon_0 = 1$.  To investigate the capability of the parameter selection strategy \eqref{eqn: MA}, we compare it with the selection technique \eqref{eqn: SQH}, which was also used in \cite{Hofmann2022}. The corresponding training results averaged over multiple, differently seeded weight initializations are reported in Table \ref{Tab: 1} and Table \ref{Tab: 2}, together with the corresponding standard deviation. Here, 'ACC' refers to the test accuracy based on $10.000$ images in the test set and 'Linesearch steps' refers to the average total number of inference steps necessary to determine a suitable augmentation parameter $\epsilon_k$ in each iteration. Further, we report the average percentage level of sparsity in the parameters of the trained models in Table \ref{Tab: 1} and  Table \ref{Tab: 2}.
For the parameter setting of the full-batch scheme, we choose $ \zeta = 0.01$, $\mu = 7$, similar to  \cite{Hofmann2022}, which can be interpreted as a coarse grid search for an appropriate augmentation parameter. On the one hand, this choice allows large updates in the learning parameters, positively affecting the algorithm's convergence on an iteration basis. On the other hand, it can result in an underestimation of $\hat{\epsilon}_k$ and thus a higher number of optimization steps discarded in the line search procedure. One can see in Table \ref{Tab: 1} and Table \ref{Tab: 2} that the latter effect is mitigated by the proposed strategy \eqref{eqn: MA}. In particular, our proposed estimation strategy for $\epsilon_k$ more than halves the required number of optimization steps compared to the standard SQH \eqref{eqn: SQH} approach by simultaneously retaining the testing accuracy and loss level achieved with \eqref{eqn: SQH}. Further, since Theorem \ref{Theorem-sufficient decrease} together with Lemma \ref{Lemma-Convergence} guarantees that the bSQH scheme imposes a decrease on the loss function method even in the non-continuous setting of the sparsity enforcing regularizer \eqref{L2L0 regularization} it can be readily applied for generating sparse CNNs. In the given case, the application of the bSQH scheme resulted in a LeNet-5 which uses only about $15 \%$ of its parameters for classification by simultaneously achieving a level of Accuracy on the test set, which is similar to that of a network using about $60 \%$ of its parameters.

\begin{table}[H]
\footnotesize
\caption{Training results for the bSQH method for classification based on the MNIST dataset and Elastic-Net regularization.}\label{Tab: 1}
\begin{center}%
\begin{tabular}{cccccc}
\hline
$\eqref{eqn: MA}$ & ACC         & $J_B$  & Sp($\%$) & Lines. steps \\ \hline
$\omega = 3$   & $98.78 \pm 0.04 $ & $0.0462 \pm 0.0023 $   & $39.79  \pm 4.77$  & $1536 \pm 29$                \\ 
$\omega = 5$  & $98.80 \pm 0.07 $ & $0.0471 \pm 0.0018 $    & $33.04 \pm 6.70$ & $1587 \pm 4$             \\ 
$\omega = 7$ & $98.80 \pm 0.04 $ & $0.0498 \pm 0.0021 $    & $29.72 \pm 5.96$  & $1615 \pm 15$                  \\ \hline
\eqref{eqn: SQH}   & $98.82 \pm 0.07 $ & $0.0444 \pm 0.0015 $    & $36.73 \pm 5.43$ & $3547 \pm 1$              \\ \hline
\end{tabular}
\end{center}
\end{table}

\begin{table}[H]
\footnotesize
\caption{Training results for the bSQH method for classification based on the MNIST dataset using the $L^0$-based regularization \eqref{L2L0 regularization}.}\label{Tab: 2}
\begin{center}%
\begin{tabular}{cccccc}
\hline
$\eqref{eqn: MA}$ & ACC         & $J_B$  & Sp($\%$) & Lines. steps \\ \hline
$\omega = 3$   & $98.67 \pm 0.02 $ & $0.1752 \pm 0.0242 $   & $82.07  \pm 2.73$  & $1573 \pm 6$                \\ 
$\omega = 5$  & $98.75 \pm 0.04 $ & $0.1726 \pm 0.0161 $    & $82.32 \pm 1.82$ & $1617 \pm 13$             \\ 
$\omega = 7$ & $98.74 \pm 0.05 $ & $0.2005 \pm 0.0224 $    & $79.16 \pm 2.51$  & $1636 \pm 14$                  \\ \hline
\eqref{eqn: SQH}   & $98.75 \pm 0.05 $ & $0.1529 \pm 0.0512 $    & $84.54 \pm 2.01$ & $3553 \pm 3$              \\ \hline
\end{tabular}
\end{center}
\end{table}

\subsection{ Mini-batch bSQH training }

\begin{figure}[H]
    \centering
    \includegraphics[width=0.9\linewidth]{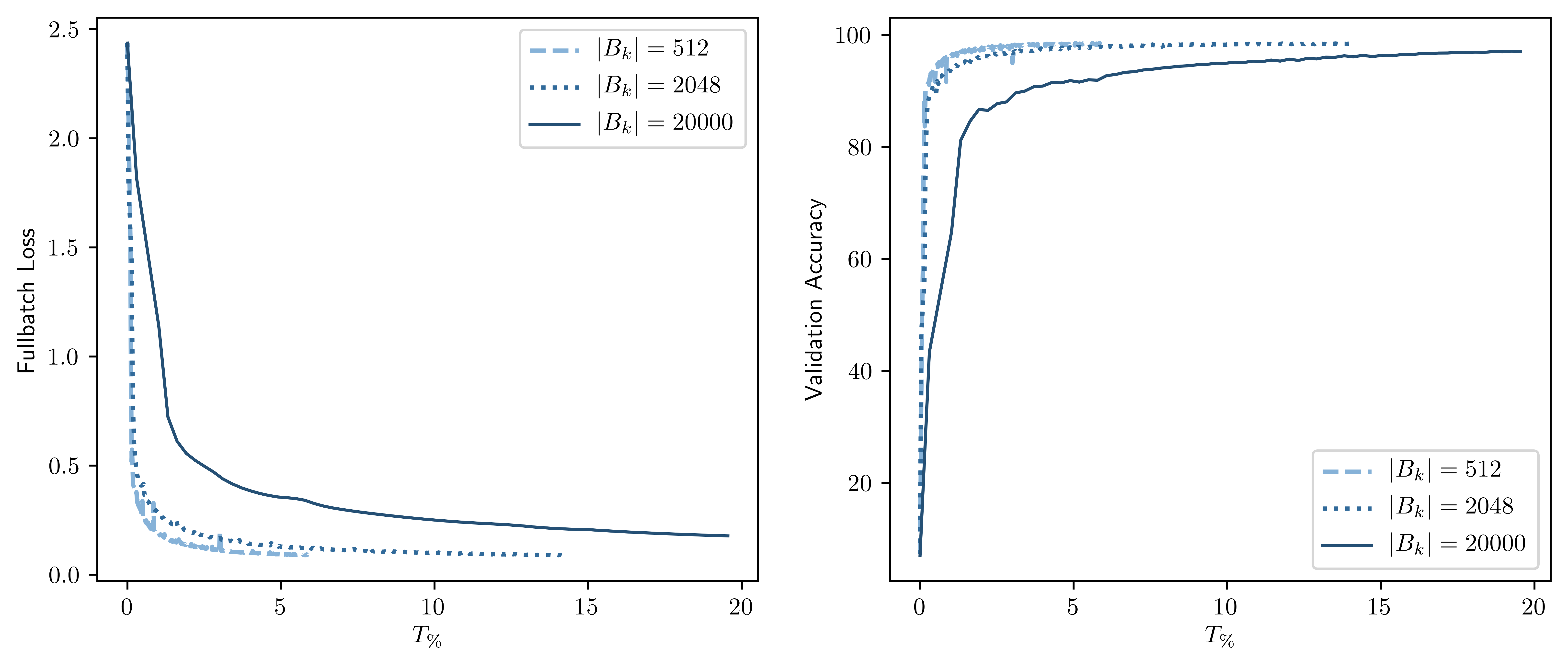}
\caption{ Loss and accuracy curves for the first 1000 iterations and with respect to the relative GPU-time of the bSQH method with $\mu = 1.1$ and the strategy \eqref{eqn: MA} with $\omega = 5$ and $\zeta = 1$. }
    \label{fig:2}
\end{figure}

In the following, the mini-batch variant of Algorithm \ref{algoSQH} is investigated. For this purpose the previously defined Elastic-Net regularized problem is considered and the standard LeNet-5 is trained for a maximum of $k_{max} = 2000$ iterations and different mini-batch sizes, that is, $ |B_k|$ $= 512$, $\;1024$, $\;2048$, $\;20.000$, where the latter corresponds to the full-batch mode of Algorithm \ref{algoSQH}. For the hyperparameters of the bSQH scheme, we set $\mu = 1.1$ and determine estimates for the augmentation parameter at each iteration using \eqref{eqn: MA}, where $\omega = 5$ and  $\zeta = 1$. This choice is motivated by the observation that a reduction of the variance in $\epsilon_k$ leads to a reduced variance in the mini-batch augmented HP functions, which carries over to the updating step of the network parameters. Further, the resulting damped variation of the augmentation parameter also prevents a constant decrease of the augmentation weight over the course of training, which would otherwise lead to overfitting to the individual mini-batches at each iteration.
The learning curves are compared based on the relative GPU time required for training, that is, $T_{\%}$ denotes the percentage of training time required relative to the most time-consuming training process, which in the given case is the full-batch setting. In Figure \ref{fig:2} and for small mini-batch sizes, the variance induced by the mini-batch approach is apparent. However, when increasing the mini-batch size, due to Lemma \ref{Lemma- Unbiased Hamiltonian estimator} and the law of large numbers, the noise in the full-batch loss and validation accuracy is reduced, resulting in an almost monotonically decreasing sequence of cost functional values, similar to the full-batch setting. This provides evidence that the minimizing properties of the full-batch bSQH algorithm corresponding to Lemma \ref{Lemma-Convergence} hold at least approximately in the mini-batch setting. Based on our results, the mini-batch approach of bSQH has favourable convergence properties with respect to the GPU time and thus yields an appropriate learning strategy for training CNNs in a stochastic framework.

It is apparent that the Assumptions of Theorem \ref{Theorem-Convergence bSQH} can easily be satisfied in the full-batch framework when considering the training of a reformulated network with layer functions given by \eqref{Control affine MLP}. For investigating these properties numerically in the mini-batch bSQH setting, we define 
$$
\Delta h : = \sum_{l = 0}^{L-1} \max_{w \in U_l} \; \cH^B_{l,\eps_k}(\By_l^{(k)}, w, u^{(k)}_l)  - \cH^B_{l,\eps_k}(\By_l^{(k)},u^{(k+1)}_l, u^{(k)}_l),\quad \Delta \Bu : = \sum_{i = k - 5}^{k}\normiii{\Bu^{(i+1)} - \Bu^{(i)}},
$$
with $k = 2000$. While the value of $\Delta h$ serves as an estimate for \eqref{Ch 4.2: almost optimal iterate}, we have that $\Delta \Bu$ enables investigating the convergence properties of the parameter sequence $(\Bu^{(k)})$ generated by the bSQH scheme. For our mini-batch framework with fixed mini-batch sizes, we observe in Table \ref{tab:3} that the conditions for Theorem \ref{Theorem-Convergence bSQH} are satisfied at least approximately. Further, contributed to Lemma \ref{Lemma- Unbiased Hamiltonian estimator}, we observe that these convergence properties are improved if the mini-batch size is increased. 
\begin{table}[H]
\footnotesize
\caption{Estimates for the validation of the assumption stated in Theorem \ref{Theorem-Convergence bSQH} when considering the the mini-batch bSQH framework.}\label{tab:3}
\begin{center}%
\begin{tabular}{cccc}
\hline
& $|B_k| = 512$  & $|B_k| =1024$  & $|B_k| =2048$ \\ \hline
$\Delta h_l$   & $1.49\mathrm{e}{-03} \pm 4.60\mathrm{e}{-04}$   & $8.46\mathrm{e}{-04} \pm 2.14\mathrm{e}{-04}$ & $9.28\mathrm{e}{-05} \pm 1.08\mathrm{e}{-04}$ \\ \hline
$\Delta \Bu$  & $6.02\mathrm{e}{-04} \pm 2.30\mathrm{e}{-04}$ & $4.80\mathrm{e}{-04} \pm 2.02\mathrm{e}{-04}$    & $1.64\mathrm{e}{-04} \pm 3.07\mathrm{e}{-05}$            \\ \hline
\end{tabular}
\end{center}
\end{table}
\subsection{ The bSQH method for classification of medical images }
\begin{figure}[tbhp]
    \centering
    \includegraphics[width=0.5\linewidth]{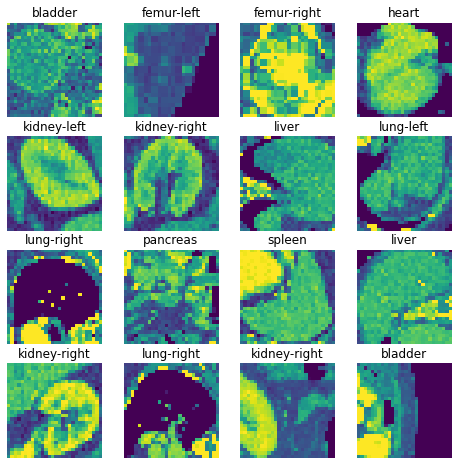}
    \caption{Samples from the OrganAMNIST dataset with corresponding labels.}
    \label{fig: Samples OrganMNIST}
\end{figure}

In this Section, we apply the mini-batch version of the bSQH method to the OrganAMNIST dataset, which is provided as part of the medMNIST database \cite{Yang2021medMNIST, Yang2023medMNIST} and consists of a collection of computed tomography (CT) images of organs; see Figure \ref{fig: Samples OrganMNIST} for samples with corresponding labels. The dataset consists of $58.850$ images, which are split into eleven categories. For solving the multiclassification task, the dataset is split into $|B| = 34.581$ samples for training, $6.491$ images for validation and $17.778$ for testing.

The CNN architecture we consider is a modification of the LeNet-5 architecture; see Table \ref{tab:modified LeNet}. These modifications involve an increase in the number of channels for each of the two convolutional layers, to allow for the extraction of more complex features, and the addition of trainable batch normalization operations (BN), aiding to reduce the internal covariate shift between layers. Each average-pooling operation is exchanged with a max-pooling operation to capture local feature information better. Additionally, we exchange the $\tanh$ activation with the widely used ReLU function, which is much simpler computationally and introduces sparsity into the feature maps by zeroing out negative activations, leading to faster forward and backward passes through the network and reduced memory requirement.

\begin{table}[tbhp]
\caption{Modified LeNet architecture description.}\label{tab:modified LeNet}
\footnotesize
\begin{center}
\begin{tabular}{ccccccc}
\hline
\multicolumn{2}{c}{Layer}           & Channels & Size  & Kernel size & Stride / Padding & $\sigma$ \\ \hline
Input   & -                         & 1        & 28x28 & -           & -                & -          \\ \hline
$\cW_1$     & Conv.                     & 16       & 28x28 & 5x5         & 2/2              & ReLU       \\
$\cP_1$ & BN + Max.Pooling & 16       & 14x14 & 2x2         & 2/-              & -          \\ \hline
$\cW_2$     & Conv.                     & 32       & 10x10 & 5x5         & 2/-              & ReLU       \\
$\cP_2$ & BN. + Max.Pooling & 32       & 5x5   & 2x2         & 2/-              & -          \\ \hline
$\cW_3$     & FC                        & -        & 800   & -           & -                & ReLU       \\
$\cP_3$ & BN               & -        & -     & -           & -                & -          \\ \hline
$\cW_4$     & FC                        & -        & 400   & -           & -                & ReLU       \\
$\cP_4$ & BN               & -        & -     & -           & -                & -          \\ \hline
$\cW_5$     & FC                        & -        & 200   & -           & -                & -          \\
Output  & -                         & -        & 11    & -           & -                & Softmax  \\ \hline 
\end{tabular}%
\end{center}
\end{table}

\begin{figure}[H]
    \centering
    \includegraphics[width=0.75\linewidth]{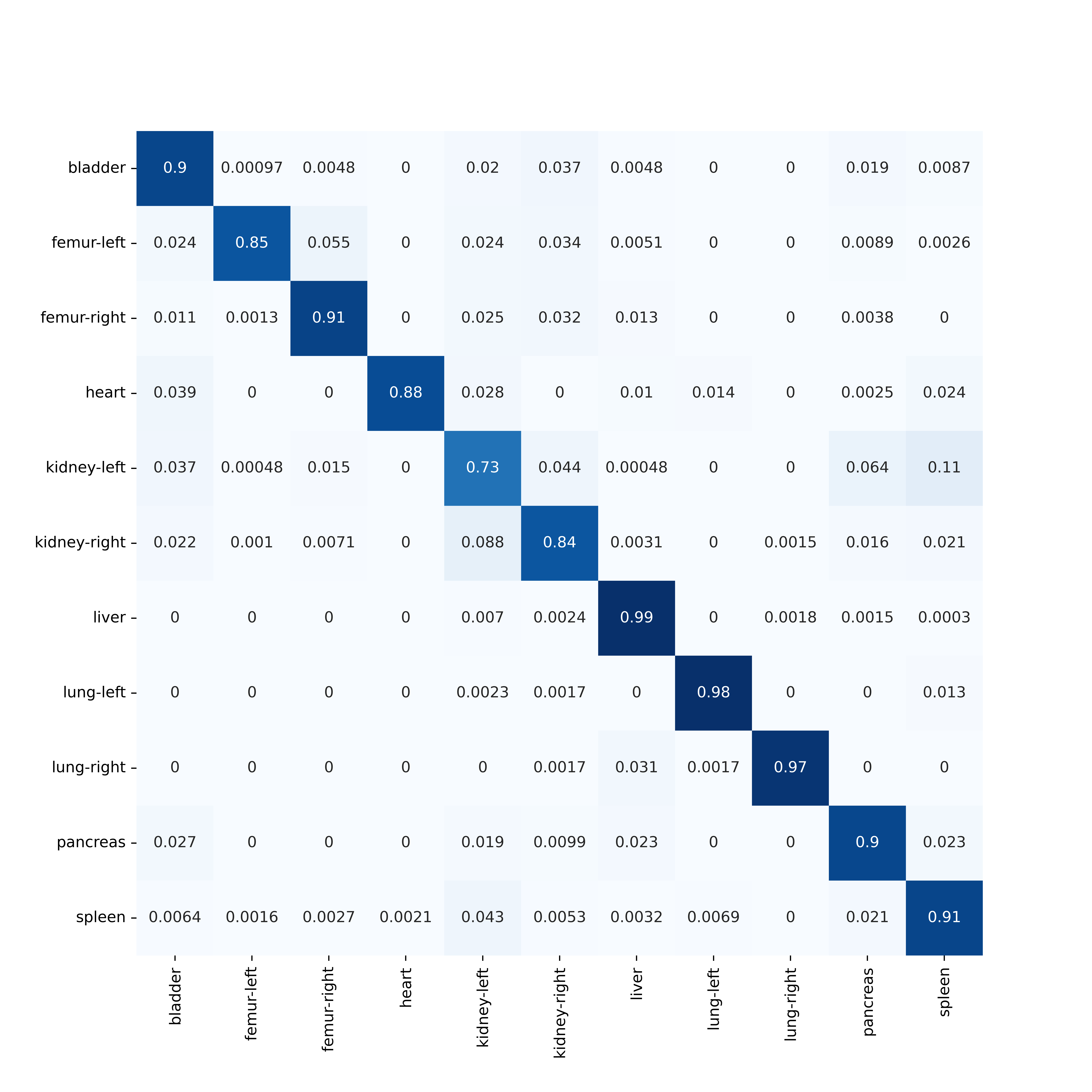}
    \caption{Confusion matrix for the modified LeNet classifier trained with mini-batch bSQH and a batch size $|B_k| = 64$. }
    \label{fig: Confusion matrix medMNSIT}
\end{figure}

The empirical risk minimization problem is formulated by combining a weighted cross entropy loss function with the sparsity enforcing regularizer \eqref{L2L0 regularization} where $\alpha = 0.99$ and a regularization weight  $\rho = 7.5\mathrm{e}{-03} $ is chosen. The weighting of the loss function is chosen to compensate for a moderate data imbalance within the training set. That is for the class $i$ we choose a weight as follows $
w_i = \exp\left( \frac{|B|}{11 \cdot c_i} \right)$, where $c_i$ denotes the total number of samples corresponding to the $i$-th class. The mini-batch bSQH method is applied with a maximum number of $k_{\max} = 8.000$ iterations, using $\mu = 1.1$ and the strategy \eqref{eqn: MA} with $\omega = 7, \zeta = 1$ together with different mini-batch sizes $|B_k| = 64,\; 128,\; 256$. We report in Table \ref{tab:medMNIST training results} the maximum testing accuracy (ACC) achieved throughout training, together with the corresponding values of the area under the receiver operating characteristic curve (AUC). Each result is the average over multiple weight initialization and is always reported together with the corresponding standard deviations. The AUC, in contrast to the ACC, is a threshold-free metric and thus less sensitive to class imbalance. Its values range from $0.5$, representing a random classifier, to $1.0$, which indicates a perfect classifier. In addition to the above evaluation metrics, we report the average percentage sparsity of the trained parameters corresponding to the reported ACC and AUC values. The confusion matrix for the best classifier results corresponding to Table \ref{tab:medMNIST training results} is given in Figure \ref{fig: Confusion matrix medMNSIT}. The confusion matrix reveals that, in particular, the classes corresponding to CT images of the kidney cause a decrease in the network's accuracy. This is due to strong similarities in the images of 'kidney-left' and 'kidney-right', possibly calling for a more elaborate network architecture. Nevertheless, and in particular, for a LeNet-5 type network architecture, the achieved results are satisfactory and can keep up with the benchmark results reported in \cite{Yang2023medMNIST}. In this context, it is important to emphasize that these results come from a model with a high level of sparsity that has only about $48.000$ trainable parameters. In contrast, models such as ResNet-18 have over a million trainable parameters and only achieve a $3\%$ increase in accuracy on the given dataset.

In addition, our results support our claim that the presented work can be extended to a framework, where instead of differentiability in Assumption \ref{A2} only Lipschitz continuity is considered. This opens up the path for an application of the bSQH scheme to state-of-the-art network architectures utilizing ReLU activations and max-pooling operations.  
\begin{table}[tbhp]
\footnotesize
\caption{Training results of the mini-batch bSQH method for classification of CT-images in the OrganAMNIST dataset}\label{tab:medMNIST training results}
\begin{center}%
\begin{tabular}{ccccc}
\hline
$B_k$ & ACC (max.)        & AUC (max.)  & Sp($\%$) \\ \hline
64   & $89.88 \pm 0.35 \,(90.34) $ & $0.9923 \pm 0.0014 \,(0.9933)$   & $78.63 \pm 3.27$                \\ \hline
128   & $89.62 \pm 0.53 \,(90.16)$ & $0.9928 \pm 0.0020 \,(0.9940)$    & $75.62 \pm 4.52$             \\ \hline
256   & $89.12 \pm 0.74 \,(89.99)$ & $0.9911 \pm 0.0034 \,(0.9942)$    & $71.48 \pm 10.90$              \\ \hline
\end{tabular}
\end{center}
\end{table}
\newpage
	\section*{Conclusion}

A supervised learning method for training convolutional neural networks (CNNs) based on a discrete-time Pontryagin maximum principle (PMP) was analyzed theoretically and validated numerically. The corresponding batch sequential quadratic Hamiltonian (bSQH) algorithm
consists of the iterative application of forward and backward sweeps with layerwise approximate maximization of an augmented Hamiltonian function. Using a reformulation of the underlying CNN, the updating step with the augmented Hamiltonian was made explicit, which enabled the application of our algorithm to $L^0$-regularized supervised learning problems. 

Convergence results for the bSQH scheme were presented, and the effectiveness of its mini-batch variant was investigated numerically. Further, it was shown that the proposed bSQH algorithm in its full-batch and mini-batch modes has superior computational performance compared to a similar existing PMP-based scheme. Results of numerical experiments demonstrated the ability of the bSQH algorithm for sparse training of CNNs in image classification tasks by including a sparsity-enforcing $L^0$-regularization.

\newpage

\bibliographystyle{siamplain}
\bibliography{Literature}

\end{document}